\documentclass[12pt]{amsart}
\usepackage{amscd,amsxtra,amssymb}

\newtheorem{theorem}{Theorem}[section]
\newtheorem{corollary}[theorem]{Corollary}
\newtheorem{lemma}[theorem]{Lemma}
\newtheorem{proposition}[theorem]{Proposition}

\theoremstyle{definition}
\newtheorem{definition}[theorem]{Definition}
\newtheorem{remark}[theorem]{Remark}
\newtheorem{example}[theorem]{Example}

\DeclareMathOperator{\Dbcoh}{\mathsf{D^{b}_{coh}}}
\DeclareMathOperator{\FM}{\mathsf{FM}}
\DeclareMathOperator{\SL}{\mathsf{SL}}
\DeclareMathOperator{\rk}{rk}
\DeclareMathOperator{\coker}{coker}
\DeclareMathOperator{\tors}{\mathsf{tors}}
\DeclareMathOperator{\Coh}{\mathsf{Coh}}
\DeclareMathOperator{\Pic}{Pic}

\DeclareMathOperator{\Hom}{\mathsf{Hom}}
\DeclareMathOperator{\Ext}{\mathsf{Ext}}
\DeclareMathOperator{\GL}{\mathsf{GL}}
\DeclareMathOperator{\Aut}{\mathsf{Aut}}
\DeclareMathOperator{\End}{\mathsf{End}}
\DeclareMathOperator{\Spec}{\mathsf{Spec}}
\DeclareMathOperator{\Stab}{\mathsf{Stab}}
\newcommand{\dtens}{\stackrel{\boldsymbol{L}}{\otimes}} 
\newcommand{\ARG}{\,\cdot\,}
\input xy
\xyoption{all}

\begin{document}
\title
[Derived categories of irreducible projective curves]
{Derived categories of irreducible projective curves of arithmetic genus one
}

\author{Igor Burban}
\address{%
Johannes-Gutenberg Universit\"at Mainz,
Fachbereich Physik, Mathematik und Informatik,
Institut f\"ur Mathematik,
55099 Mainz, Germany
}
\email{burban@mathematik.uni-mainz.de}

\author{Bernd Kreu{\ss}ler}
\address{%
Mary Immaculate College, South Circular Road, Limerick, Ireland
}
\email{bernd.kreussler@mic.ul.ie}


\date{}

\dedicatory{Dedicated to Yuriy Drozd on his sixtieth birthday} 

\begin{center}
{\large
}
\end{center}

\begin{abstract}
We investigate the bounded derived category of coherent sheaves on irreducible
singular projective curves of arithmetic genus one.  
A description of the group of exact auto-equivalences and the set of all
t-structures of this category is given.
We describe the moduli space of stability conditions, obtain a complete
classification of all spherical objects in this category
and show that the group of exact auto-equivalences acts transitively on
them. 
Harder-Narasimhan filtrations in the sense of Bridgeland are used as our main
technical tool.
\end{abstract}


\maketitle

\section{Introduction}

The purpose of this paper is to study the structure of the bounded derived
category $\Dbcoh(\boldsymbol{E})$ of coherent sheaves on a singular irreducible
projective curve $\boldsymbol{E}$ of arithmetic genus one. 

In the smooth case, such structure results are easily obtained from Atiyah's
description \cite{Atiyah} of indecomposable vector bundles over elliptic
curves. However, if $\boldsymbol{E}$ has a node or a cusp, some crucial 
properties fail to hold. This is illustrated by the following table. 
\begin{center}
  \begin{tabular}[t]{p{6cm}|c|c}
     &smooth&singular\\ \hline
    homological dimension of $\Coh_{\boldsymbol{E}}$
     &$1$&$\infty$\\ \hline
    Serre duality holds&in general&\multicolumn{1}{p{3cm}}{
    with one object being perfect}\\ \hline
    torsion free implies locally free&yes&no\\ \hline
    indecomposable coherent sheaves are semi-stable&yes&no\\ \hline
    any indecomposable complex is isomorphic to a shift of a sheaf&yes&no\\
    \hline 
  \end{tabular}
\end{center}
Despite these difficulties, the main goal of this article is to find the common
features between the smooth and the singular case.  A list of such can be
found in Remark \ref{rem:common}.

In Section \ref{sec:background}, we review the smooth case and highlight where
the properties mentioned above are used. Our approach was inspired by
\cite{LenzingMeltzer}.  
Atiyah's algorithm to construct indecomposable vector bundles of any slope can
be understood as an application of a sequence of twist functors with spherical
objects. From this point of view, Atiyah shows that any indecomposable object
of $\Dbcoh(\boldsymbol{E})$ is the image of an indecomposable torsion sheaf
under an exact auto-equivalence of $\Dbcoh(\boldsymbol{E})$.

In the case of a singular Weierstra{\ss} curve $\boldsymbol{E}$, as our main
technical tool we use Harder-Narasimhan filtrations in
$\Dbcoh(\boldsymbol{E})$, which were introduced by Bridgeland
\cite{Stability}. Their general properties are studied in Section
\ref{sec:HNF}.  

The key result of Section \ref{sec:dercat} is the preservation of stability
under Seidel-Thomas twists \cite{SeidelThomas} with spherical objects.  This
allows us to show that, like in the smooth case, any category of semi-stable
sheaves with fixed slope is equivalent to the category of coherent torsion
sheaves on $\boldsymbol{E}$. 

In the case of slope zero, this was shown in our previous work
\cite{BurbanKreussler}. For the nodal case, an explicit description of
semi-stable sheaves of degree zero via \'etale coverings was given there as
well.   
A combinatorial description of semi-stable sheaves of arbitrary slope over
a nodal cubic curve was found by Mozgovoy \cite{Mozgovoy}.
On the other hand, a classification of all indecomposable objects of
$\Dbcoh(\boldsymbol{E})$ was presented in \cite{BurbanDrozd}. A description of
the Harder-Narasimhan filtrations in terms of this classification is a task for
future work. 

However, if the singular point of $\boldsymbol{E}$ is a cusp, the description
of all indecomposable coherent torsion sheaves is a wild problem in the sense
of representation theory, see for example \cite{Drozd72}. Nevertheless,
stable vector bundles on a cuspidal cubic have been classified by Bodnarchuk
and Drozd \cite{Lesya}. 

It turns out that semi-stable sheaves of infinite homological dimension are
particularly important, because only such sheaves appear as Harder-Narasimhan
factors of indecomposable objects in $\Dbcoh(\boldsymbol{E})$ which are not
semi-stable (Proposition \ref{prop:extreme}). 

The main result (Proposition \ref{prop:spherical}) of Section \ref{sec:dercat}
is the answer to a question of Polishchuk, who asked in \cite{YangBaxter},
Section 1.4, for a description of all spherical objects on $\boldsymbol{E}$. 
We also prove that the group of exact auto-equivalences of
$\Dbcoh(\boldsymbol{E})$ acts transitively on the set of spherical objects.

In Section \ref{sec:tstruc} we study $t$-structures on $\Dbcoh(\boldsymbol{E})$
and stability conditions in the sense of \cite{Stability}. 
We completely classify all $t$-structures on this category (Theorem
\ref{thm:tstruc}). This allows us to 
deduce a description of the group of exact auto-equivalences of
$\Dbcoh(\boldsymbol{E})$ (Corollary \ref{cor:auto}).
As a second application, we calculate Bridgeland's moduli space of
stability conditions on $\boldsymbol{E}$ (Proposition \ref{prop:stabmod}). 

The hearts $\mathsf{D}(\theta,\theta+1)$ of the $t$-structures constructed in
Section \ref{sec:tstruc} are finite-dimensional non-Noetherian Abelian
categories of infinite global dimension.  
In the case of a smooth elliptic curve, this category is equivalent to the
category  of holomorphic vector bundles on a non-commutative torus in the
sense of Polishchuk and Schwarz \cite{PolSchw}, Proposition 3.9.  
It is an interesting problem to find such a differential-geometric 
interpretation of these Abelian categories in the case of singular
Weierstra{\ss} curves. 

Using the technique of Harder-Narasimhan filtrations, we gain new insight into
the classification of indecomposable complexes, which was obtained in
\cite{BurbanDrozd}. 
It seems plausible that similar methods can be applied to study the derived
category of representations of certain derived tame associative algebras, such
as gentle algebras, skew-gentle algebras or degenerated tubular algebras, see
for example \cite{BuDro}.
The study of Harder-Narasimhan filtrations in conjunction with the action of
the group of exact auto-equivalences of the derived category should provide new
insight into the combinatorics of indecomposable objects in these derived
categories. 

\textbf{Notation.} We fix an algebraically closed field $\boldsymbol{k}$ of
characteristic zero. By $\boldsymbol{E}$ we always denote a Weierstra{\ss}
curve. This is a reduced irreducible curve of arithmetic genus one, isomorphic
to a cubic curve in the projective plane. If not smooth, it has precisely one
singular point $s\in\boldsymbol{E}$, which can be a node or a cusp.
If $x\in\boldsymbol{E}$ is arbitrary, we denote by $\boldsymbol{k}(x)$ the
residue field of $x$ and consider it as a sky-scraper sheaf supported at $x$.
By $\Dbcoh(\boldsymbol{E})$ we denote the derived category of complexes of
$\mathcal{O}_{\boldsymbol{E}}$-modules whose cohomology sheaves are coherent
and which are non-zero only in finitely many degrees.

\textbf{Acknowledgement.} The first-named author would like to thank
Max-Planck-Institut f\"ur Mathematik in Bonn for financial support. 
Both authors would like to thank Yuriy Drozd, Daniel Huybrechts, Bernhard
Keller, Rapha\"el Rouquier and  Olivier  Schiffmann  for helpful discussions, 
and the referee for his or her constructive comments.

\section{Background: the smooth case}\label{sec:background}

The purpose of this section is to recall well-known results about the structure
of the bounded derived category of coherent sheaves over a smooth elliptic
curve. Proofs of most of these results can be found in \cite{Atiyah},
\cite{Oda}, \cite{LenzingMeltzer} and \cite{Tu}. 
The focus of our presentation is on the features and techniques
which are essential in the singular case as well. 
At the end of this section we highlight the main differences between the smooth
and the singular case. It becomes clear that the failure of Serre duality is
the main reason why the proofs and even the formulation of some of the main
results do not carry over to the singular case. 
The aim of the subsequent sections will then be to overcome these difficulties,
to find correct formulations which generalise to the singular case and to
highlight the common features of the bounded derived category in the smooth and
singular case. 

With the exception of subsection \ref{subsec:diff}, throughout this section
$\boldsymbol{E}$ denotes a smooth elliptic curve over $\boldsymbol{k}$.

\subsection{Homological dimension}

For any two coherent sheaves $\mathcal{F}, \mathcal{G}$ on $\boldsymbol{E}$,
Serre duality provides an isomorphism $$\Ext^{\nu}(\mathcal{F},\mathcal{G})
\cong \Ext^{1-\nu}(\mathcal{G},\mathcal{F})^{\ast}.$$ 
This follows from the usual formulation of Serre duality and the fact that any
coherent sheaf has a finite locally free resolution. As a consequence,
$\Ext^{\nu}(\mathcal{F},\mathcal{G})=0$ for any $\nu \ge 2$, which means that
$\Coh_{\boldsymbol{E}}$ has homological dimension one. This implies that any
object $X\in\Dbcoh(\boldsymbol{E})$ splits into the direct sum of appropriate
shifts of its cohomology sheaves. To see this, start with a complex
$X=(\mathcal{F}^{-1} \stackrel{f}{\longrightarrow} \mathcal{F}^{0})$
and consider the distinguished triangle in $\Dbcoh(\boldsymbol{E})$ 
$$\ker(f)[1] \rightarrow X \rightarrow \coker(f) \stackrel{\xi}{\rightarrow}
\ker(f)[2].$$ 
Because $\xi\in\Hom(\coker(f),\ker(f)[2]) = \Ext^{2}(\coker(f), \ker(f)) =0$,
we obtain $X\cong \ker(f)[1] \oplus \coker(f)$. Using the same idea we can
proceed by induction to get the claim.

\subsection{Indecomposable sheaves are semi-stable}

It is well-known that any coherent sheaf $\mathcal{F}\in\Coh_{\boldsymbol{E}}$
has a Harder-Narasimhan filtration
$$0\subset \mathcal{F}_{n} \subset \ldots \subset \mathcal{F}_{1} \subset
\mathcal{F}_{0} = \mathcal{F}$$ 
whose factors $\mathcal{A}_{\nu} := \mathcal{F}_{\nu}/\mathcal{F}_{\nu+1}$ are
semi-stable with decreasing slopes $\mu(\mathcal{A}_{n})>
\mu(\mathcal{A}_{n-1}) > \ldots > \mu(\mathcal{A}_{0})$.  
Using the definition of semi-stability, this implies
$\Hom(\mathcal{A}_{\nu+i}, \mathcal{A}_{\nu}) 
= 0$ for all $\nu\ge 0$ and $i>0$. Therefore,
$\Ext^{1}(\mathcal{A}_{0},\mathcal{F}_{1}) \cong 
\Hom(\mathcal{F}_{1}, \mathcal{A}_{0})^{\ast} =0$, and the exact sequence
$0\rightarrow \mathcal{F}_{1} \rightarrow \mathcal{F} \rightarrow
\mathcal{A}_{0} \rightarrow 0$ must split. 
In particular, if $\mathcal{F}$ is indecomposable, we have $\mathcal{F}_{1}=0$
and $\mathcal{F}\cong \mathcal{A}_{0}$ and $\mathcal{F}$  is semi-stable.

\subsection{Jordan-H\"older factors}

The full sub-category of $\Coh_{\boldsymbol{E}}$ whose objects are the
semi-stable sheaves of a fixed slope is an Abelian category in which any object
has a Jordan-H\"older filtration with stable factors. 
If $\mathcal{F}$ and $\mathcal{G}$ are non-isomorphic stable sheaves which
have the same slope, we have $\Hom(\mathcal{F},\mathcal{G})=0$. Based on this
fact, in the same way as before, we can deduce that an indecomposable
semi-stable sheaf has all its Jordan-H\"older factors isomorphic to each
other.

\subsection{Simple is stable}

It is well-known that any stable sheaf $\mathcal{F}$ is simple,
i.e. $\Hom(\mathcal{F},\mathcal{F}) \cong \boldsymbol{k}$. On a smooth
elliptic curve, the converse is true as well, which equips us with a useful
homological characterisation of stability.  

To see that simple implies stable, we suppose for a contradiction that
$\mathcal{F}$ is simple but not stable. This implies the existence of an
epimorphism $\mathcal{F}\rightarrow \mathcal{G}$ with $\mathcal{G}$ stable and
$\mu(\mathcal{F})\ge \mu(\mathcal{G})$. Serre duality implies  
$\dim \Ext^{1}(\mathcal{G},\mathcal{F}) = \dim \Hom(\mathcal{F},\mathcal{G}) >
0$, hence, $\chi(\mathcal{G},\mathcal{F}) := \dim
\Hom(\mathcal{G},\mathcal{F}) - \dim \Ext^{1}(\mathcal{G},\mathcal{F}) < \dim
\Hom(\mathcal{G},\mathcal{F})$. Riemann-Roch gives
$\chi(\mathcal{G},\mathcal{F}) = (\mu(\mathcal{F}) -
\mu(\mathcal{G}))/\rk(\mathcal{F})\rk(\mathcal{G}) > 0$, hence
$\Hom(\mathcal{G},\mathcal{F})\ne 0$. But this produces a 
non-zero composition $\mathcal{F}\rightarrow \mathcal{G} \rightarrow
\mathcal{F}$ which is not an isomorphism, in contradiction to the assumption
that $\mathcal{F}$ was simple.

\subsection{Classification}

Atiyah \cite{Atiyah} gave a description of all stable sheaves with a fixed
slope in the form $\mathcal{E}(r,d)\otimes \mathcal{L}$, where $\mathcal{L}$
is a line bundle of degree zero and $\mathcal{E}(r,d)$ is a particular stable
bundle of the fixed slope. The bundle $\mathcal{E}(r,d)$ depends on the choice
of a base point $p_{0}\in\boldsymbol{E}$ and 
its construction reflects the Euclidean algorithm on the pair $(\rk,\deg)$. 
We look at this description from a slightly different perspective. We use the
twist functors $T_{\mathcal{O}}$ and $T_{\boldsymbol{k}(p_{0})}$, which were
constructed by Seidel and Thomas \cite{SeidelThomas} (see also \cite{Meltzer}).
They act as equivalences on $\Dbcoh(\boldsymbol{E})$ and, hence, preserve
stability. 

A stable sheaf of rank $r$ and degree $d$ is sent by $T_{\mathcal{O}}$ to one
with $(\rk,\deg)$ equal to $(r-d,d)$. If $r<d$ this is a shift of a stable
sheaf. The functor $T_{\boldsymbol{k}(p_{0})}$ sends the pair $(r,d)$ to
$(r,r+d)$ and its inverse sends it to $(r,d-r)$.  Therefore, if we follow the
Euclidean algorithm, we find a composition of such functors which provides an
equivalence between the category of stable sheaves with slope $d/r$ and the
category of simple torsion sheaves. Such sheaves are precisely the structure
sheaves of closed points $\boldsymbol{k}(x)$, $x\in\boldsymbol{E}$. They are
considered to be stable with slope $\infty$.  

More generally, this procedure provides an equivalence between the category of
semi-stable sheaves of rank $r$ and degree $d$ with the category of torsion
sheaves of length equal to $\gcd(r,d)$. This shows, in particular, that the
Abelian category of semi-stable sheaves with fixed slope is equivalent to the
category of coherent torsion sheaves.

\subsection{Auto-equivalences}

By $\Aut(\Dbcoh(\boldsymbol{E}))$ we denote the group of all exact
auto-equivalences of the triangulated category $\Dbcoh(\boldsymbol{E})$. This
group acts on the Grothendieck group $\mathsf{K}(\boldsymbol{E}) \cong
\mathsf{K}(\Dbcoh(\boldsymbol{E}))$.
As the kernel of the Chern character is the radical of the Euler-form  
$\langle X,Y \rangle = \dim(\Hom(X,Y)) - \dim(\Hom(X,Y[1])$
which is invariant under this action, it induces an action on the even
cohomology  $H^{2\ast}(\boldsymbol{E}, \mathbb{Z}) \cong \mathbb{Z}^{2}$.  
Because $\dim(\Hom(\mathcal{F},\mathcal{G}))>0$ if and only if $\langle
\mathcal{F},\mathcal{G} \rangle >0$, provided $\mathcal{F}\not\cong
\mathcal{G}$ are stable sheaves, the induced action on $\mathbb{Z}^{2}$ is 
orientation preserving. So, we obtain a homomorphism of groups $\varphi:
\Aut(\Dbcoh(\boldsymbol{E})) \rightarrow \SL(2,\mathbb{Z})$, which is
surjective because $T_{\mathcal{O}}$ and $T_{\boldsymbol{k}(p_{0})}$ are
mapped to a pair of generators of $\SL(2,\mathbb{Z})$. Explicitly, if
$\mathbb{G}$ is an auto-equivalence, $\varphi(\mathbb{G})$ describes its
action on the pair $(\rk,\deg)$. 
To understand $\ker(\varphi)$, we observe that $\varphi(\mathbb{G}) =
\boldsymbol{1}$ implies that $\mathbb{G}$ sends a simple torsion sheaf
$\boldsymbol{k}(x)$ to some $\boldsymbol{k}(y)[2k]$, because indecomposability
is retained. By the same reason, $\mathbb{G}(\mathcal{O})$ is a shifted line
bundle of degree zero.  
However, $\Hom(\mathcal{L},\boldsymbol{k}(y)[l]) = 0$, if $\mathcal{L}$ is a
line bundle and $l\ne 0$. Hence, after composing $\mathbb{G}$ with a shift, it
sends all simple torsion sheaves to simple torsion sheaves, without a shift. 
Because $\boldsymbol{E}$ is smooth, we can apply a result of Orlov \cite{Orlov}
which says that any auto-equivalence $\mathbb{G}$ is a Fourier-Mukai transform
\cite{Mukai}. 
However, any such functor, which sends the sheaves $\boldsymbol{k}(x)$ to
torsion  sheaves of length one  is of the form
$\mathbb{G}(X)=f^{\ast}(\mathcal{L}\otimes X)$, 
where $f:\boldsymbol{E} \rightarrow \boldsymbol{E}$ is an automorphism and
$\mathcal{L}\in\Pic(\boldsymbol{E})$ a line bundle. Hence, $\ker(\varphi)$ is
generated by $\Aut(\boldsymbol{E}), \Pic^{0}(\boldsymbol{E})$  and even
shifts. This gives a complete description of the group
$\Aut(\Dbcoh(\boldsymbol{E}))$. A similar approach was used by Lenzing and
Meltzer to describe the group of exact auto-equivalences of tubular weighted
projective lines 
\cite{LenzingMeltzerAuto}.

\subsection{Difficulties in the singular case}\label{subsec:diff}

Let now $\boldsymbol{E}$ be an irreducible but singular curve of arithmetic
genus one. The technical cornerstones of the theory as described in this
section fail to be true in this case. More precisely:
\begin{itemize}
\item the category of coherent sheaves $\Coh_{\boldsymbol{E}}$ has infinite
  homological dimension; 
\item there exist indecomposable complexes in $\Dbcoh(\boldsymbol{E})$ which
  are not just shifted sheaves, see  \cite{BurbanDrozd}, section 3;
\item Serre duality fails to be true in general;
\item not all indecomposable vector bundles are semi-stable;
\item there exist indecomposable coherent sheaves which are neither torsion
  sheaves nor torsion free sheaves, see \cite{BurbanDrozd}.
\end{itemize}

Most of the trouble is caused by the failure of Serre duality.
The basic example is the following. Suppose, $s\in\boldsymbol{E}$ is a node,
then 
$$\Hom(\boldsymbol{k}(s), \boldsymbol{k}(s)) \cong \boldsymbol{k}\quad
\text{ and }\quad
\Ext^{1}(\boldsymbol{k}(s), \boldsymbol{k}(s)) \cong \boldsymbol{k}^{2}.$$
Serre duality is available, only if at least one of the two sheaves involved
has finite homological dimension. This might suggest that replacing
$\Dbcoh(\boldsymbol{E})$ by the sub-category of perfect complexes would solve
most of the problems. But see Remark \ref{rem:notperfect}.

In the subsequent sections we overcome these difficulties and point out the
similarities between the smooth and the singular case.

\section{Harder-Narasimhan filtrations}\label{sec:HNF}

Throughout this section, $\boldsymbol{E}$ denotes an irreducible reduced
projective curve over $\boldsymbol{k}$ of arithmetic genus one.
The notion of stability of coherent torsion free sheaves on an irreducible
curve is usually defined with the aid of the slope function
$\mu(\ARG)=\deg(\ARG)/\rk(\ARG)$. To use the phase function instead is
equivalent, but better adapted for the generalisation to derived categories
described below. 
By definition, the \emph{phase} $\varphi(\mathcal{F})$ of a non-zero
coherent sheaf $\mathcal{F}$ is the unique number which satisfies $0 <
\varphi(\mathcal{F})\le 1$ and $m(\mathcal{F}) \exp(\pi
i\varphi(\mathcal{F})) = -\deg(\mathcal{F}) + i \rk(\mathcal{F})$, where
$m(\mathcal{F})$ is a positive real number, called the \emph{mass} of the
sheaf $\mathcal{F}$. 
In particular, $\varphi(\mathcal{O}) = 1/2$ and all non-zero torsion sheaves
have phase one.
A torsion free coherent sheaf $\mathcal{F}$ is called semi-stable if for any
exact sequence of torsion free coherent sheaves 
$$0 \rightarrow \mathcal{E} \rightarrow \mathcal{F}
    \rightarrow \mathcal{G} \rightarrow 0$$
the inequality $\varphi(\mathcal{E}) \le \varphi(\mathcal{F})$, or
equivalently, $\varphi(\mathcal{F}) \le \varphi(\mathcal{G})$, holds. 
It is well-known \cite{Rudakov} that any torsion free coherent sheaf
$\mathcal{F}$ on a projective variety has a Harder-Narasimhan filtration  
$$0 \subset \mathcal{F}_{n} \subset \mathcal{F}_{n-1}  \cdots \subset
 \mathcal{F}_{1} \subset \mathcal{F}_{0}  = \mathcal{F},$$
which is uniquely characterised by the property that all factors
 $\mathcal{A}_{i} = \mathcal{F}_{i}/\mathcal{F}_{i+1}$ are semi-stable and
satisfy  
$$\varphi(\mathcal{A}_{n}) > \varphi(\mathcal{A}_{n-1}) > \cdots >
  \varphi(\mathcal{A}_{0}).$$

Originally, this concept of stability was introduced in the 1960s in order to
construct moduli spaces using geometric invariant theory. It could also be
seen as a method to understand the structure of the category of coherent
sheaves on a projective variety. 
By Simpson, the notion of stability was extended to coherent sheaves of
pure dimension.  
A very general approach was taken by Rudakov \cite{Rudakov}, who introduced
the notion of stability on Abelian categories. Under some finiteness
assumptions on the category, he shows the existence and uniqueness of a
Harder-Narasimhan filtration for any object of the category in question. As an
application of his work, the usual slope stability extends to the whole
category $\Coh_{\boldsymbol{E}}$ of coherent sheaves on $\boldsymbol{E}$. In
particular, any non-zero coherent sheaf has a Harder-Narasimhan filtration and
any non-zero coherent torsion sheaf on the curve $\boldsymbol{E}$ is
semi-stable.

Inspired by work of Douglas on $\Pi$-Stability for D-branes, see for example
\cite{Douglas}, it was shown by Bridgeland \cite{Stability} how to extend the
concept of  stability and Harder-Narasimhan filtration to the derived category
of coherent sheaves, or more generally, to a triangulated category. These new
ideas were merged with the ideas from \cite{Rudakov} in the paper
\cite{GRK}. 
We shall follow here the approach of Bridgeland \cite{Stability}. In Section 
\ref{sec:tstruc} we give a description of Bridgeland's moduli space of
stability conditions on the derived category of irreducible singular curves
of arithmetic genus one. However, throughout the present chapter we stick to
the classical notion of stability on the category of coherent sheaves and the
stability structure it induces on the triangulated category.

In order to generalise the concept of a Harder-Narasimhan filtration to the
category $\Dbcoh(\boldsymbol{E})$, Bridgeland \cite{Stability} extends the
definition of the phase of a sheaf to shifts of coherent sheaves by:
$$\varphi(\mathcal{F}[n]) := \varphi(\mathcal{F})+n,$$
where $\mathcal{F}\ne 0$ is a coherent sheaf on $\boldsymbol{E}$ and
$n\in\mathbb{Z}$. A complex which is non-zero at position $m$ only has,
according to this definition, phase in the interval $(-m,-m+1]$. If
$\mathcal{F}$ and $\mathcal{F}'$ are non-zero coherent sheaves and $a,b$
integers, we have the implication:
$$\varphi(\mathcal{F}[-a]) > \varphi(\mathcal{F}'[-b])
\quad\Rightarrow\quad a\le b.$$
For any $\varphi\in\mathbb{R}$ we denote by $\mathsf{P}(\varphi)$ the Abelian
category of shifted semi-stable sheaves with phase $\varphi$. Of course,
$0\in\mathsf{P}(\varphi)$ for all $\varphi$. 
If $\varphi\in(0,1]$, this is a full Abelian subcategory of
$\Coh_{\boldsymbol{E}}$. For any $\varphi\in\mathbb{R}$ we have
$\mathsf{P}(\varphi+n) = \mathsf{P}(\varphi)[n]$. A non-zero object of
$\Dbcoh(\boldsymbol{E})$ will be called \emph{semi-stable}, if it is an
element of one of the categories $\mathsf{P}(\varphi)$, $\varphi\in\mathbb{R}$.

Bridgeland's stability conditions \cite{Stability} involve so-called
central charges.  
In order to define the central charge of the standard stability condition, we
need a definition of degree and rank for arbitrary objects in
$\Dbcoh(\boldsymbol{E})$.

Let $K =\mathcal{O}_{\boldsymbol{E},\eta}$ be the field of rational
functions on the irreducible curve $\boldsymbol{E}$ with generic point
$\eta\in\boldsymbol{E}$. The base change $\eta:\Spec(K)\rightarrow
\boldsymbol{E}$ is flat, so that  $\eta^{\ast}(F)$, 
taken in the non-derived sense, is correctly defined for any
$F\in\Dbcoh(\boldsymbol{E})$. 
We define $\rk(F):=\chi(\eta^{\ast}(F))$, which is the
alternating sum of the dimensions of the cohomology spaces of the complex
$\eta^{\ast}(F)$ which are vector spaces over $K$.

In order to define the degree, we use the functor 
$$\boldsymbol{R}\Hom(\mathcal{O}_{\boldsymbol{E}},\ARG): \Dbcoh(\boldsymbol{E})
\rightarrow \Dbcoh(\boldsymbol{k}),$$
and set $\deg(F):= \chi(\boldsymbol{R}\Hom(\mathcal{O}_{\boldsymbol{E}},F))$.
Here, we denoted by $\Dbcoh(\boldsymbol{k})$ the bounded derived category of
finite dimensional vector spaces over $\boldsymbol{k}$.
For coherent sheaves, these definitions coincide with the usual
definitions of rank and degree. In particular, a torsion sheaf of length $m$
which is supported at a single point of $\boldsymbol{E}$ has rank $0$ and
degree $m$. 

These definitions imply that rank and degree are additive on distinguished
triangles in $\Dbcoh(\boldsymbol{E})$. Hence, they induce homomorphisms on the
Grothendieck group $\mathsf{K}(\Dbcoh(\boldsymbol{E}))$ of the triangulated
category $\Dbcoh(\boldsymbol{E})$, which is by definition the quotient of the
free Abelian group generated by the objects of $\Dbcoh(\boldsymbol{E})$ modulo
expressions coming from distinguished triangles. 
Recall that $\mathsf{K}_{0}(\Coh(\boldsymbol{E})) \cong
\mathsf{K}(\Dbcoh(\boldsymbol{E}))$, 
see \cite{Groth}. We denote this group by $\mathsf{K}(\boldsymbol{E})$

\begin{lemma}\label{lem:GrothGrp}
  If $\boldsymbol{E}$ is an irreducible singular curve of arithmetic genus
  one, we have $\mathsf{K}(\boldsymbol{E}) \cong \mathbb{Z}^{2}$ with
  generators $[\boldsymbol{k}(x)]$ and $[\mathcal{O}_{\boldsymbol{E}}]$. 
\end{lemma}

\begin{proof}
  Recall that the Grothendieck-Riemann-Roch Theorem, see \cite{BFM} or
  \cite{Fulton}, provides a homomorphism
  $$\tau_{\boldsymbol{E}}:\mathsf{K}(\boldsymbol{E}) \rightarrow
  A_{\ast}(\boldsymbol{E})\otimes \mathbb{Q},$$
  which depends functorially on $\boldsymbol{E}$ with respect to proper direct
  images.
  Moreover,
  $(\tau_{\boldsymbol{E}})_{\mathbb{Q}}:\mathsf{K}(\boldsymbol{E})\otimes
  \mathbb{Q} \rightarrow A_{\ast}(\boldsymbol{E})\otimes \mathbb{Q}$
  is an isomorphism, see \cite{Fulton}, Cor.\/ 18.3.2.

  If $\boldsymbol{E}$ is an irreducible singular projective curve of
  arithmetic genus one, we easily see that the Chow group
  $A_{\ast}(\boldsymbol{E})$ is isomorphic to $\mathbb{Z}^{2}$.
  The two generators are $[x]\in A_{0}(\boldsymbol{E})$ with
  $x\in\boldsymbol{E}$ and $[\boldsymbol{E}]\in A_{1}(\boldsymbol{E})$.
  Note that $[x]=[y]\in A_{0}(\boldsymbol{E})$ for any two closed points
  $x,y\in\boldsymbol{E}$, because the normalisation of $\boldsymbol{E}$ is
  $\mathbb{P}^{1}$. 
  Using \cite{Fulton}, Thm.\/ 18.3 (5), we obtain
  $\tau_{\boldsymbol{E}}(\boldsymbol{k}(x))= [x] \in A_{0}(\boldsymbol{E})$
  for any $x\in\boldsymbol{E}$. On the other hand, from \cite{Fulton}, Expl.\/
  18.3.4 (a), we obtain $\tau_{\boldsymbol{E}}(\mathcal{O}_{\boldsymbol{E}}) =
  [\boldsymbol{E}] \in A_{1}(\boldsymbol{E})$.
  Therefore, the classes of $\boldsymbol{k}(x)$ and
  $\mathcal{O}_{\boldsymbol{E}}$ define a basis of
  $\mathsf{K}(\boldsymbol{E})\otimes \mathbb{Q}$.
  However, these two classes generate the group $\mathsf{K}(\boldsymbol{E})$,
  so that it must be a free Abelian group. 
\end{proof}

The \emph{central charge} of the standard stability structure on
$\Dbcoh(\boldsymbol{E})$ is the homomorphism of Abelian groups
$$
Z: \mathsf{K}(\boldsymbol{E}) \rightarrow \mathbb{Z}\oplus i\mathbb{Z}
\subset \mathbb{C},
$$
which is given by 
$$
Z(F) := -\deg(F) + i \rk(F).
$$
If $F$ is a non-zero coherent sheaf, $Z(F)$ is a point on the ray from the
origin through $\exp(\pi i \varphi(F))$ in $\mathbb{C}$. Its distance from the
origin was called the mass of $F$. 

Although the phase $\varphi(F)$ is defined for sheaves and their shifts only,
we are able to define the slope $\mu(F)$ for any object in
$\Dbcoh(\boldsymbol{E})$ which is not equal to zero in the Grothendieck
group. Namely, the usual definition $\mu(F):=\deg(F)/\rk(F)$ gives us now a
mapping  $$\mu:\mathsf{K}(\boldsymbol{E})\setminus\{0\} \rightarrow
\mathbb{Q} \cup \{\infty\},$$ which extends the usual definition of the slope
of a sheaf. 
Because $Z(\mathcal{O}_{\boldsymbol{E}})=i$ and $Z(\boldsymbol{k}(x))=-1$,
Lemma \ref{lem:GrothGrp} implies that $Z$ is injective. Therefore, $\mu$ is
defined for any non-zero element of the Grothendieck group.

For arbitrary objects $X\in\Dbcoh(\boldsymbol{E})$ we have $Z(X[1]) = -Z(X)$,
hence $\mu(X[1]) = \mu(X)$ when defined.
In case of shifted sheaves, in contrast to the slope $\mu$, the
phase $\varphi$ keeps track of the position of this sheaf in the complex.
As an illustration, we include an example of an indecomposable object in 
$\Dbcoh(\boldsymbol{E})$ which has a zero image in the Grothendieck group. 

\begin{example}
  Let $s\in\boldsymbol{E}$ be the singular point and denote, as usual, by
  $\boldsymbol{k}(s)$ the torsion sheaf of length one which is supported at
  $s$. This sheaf does not have finite homological dimension. To see this, we
  observe first that $\Ext^{k}(\boldsymbol{k}(s), \boldsymbol{k}(s)) \cong
  H^{0}(\mathcal{E}xt^{k}(\boldsymbol{k}(s), \boldsymbol{k}(s)))$. Moreover, 
  as an 
  $\mathcal{O}_{\boldsymbol{E},s}$-module, 
  $\boldsymbol{k}(s)$ has an infinite periodic locally free resolution of the
  form 
  $$
    \cdots \stackrel{A}{\longrightarrow} \mathcal{O}_{\boldsymbol{E},s}^{2} 
           \stackrel{B}{\longrightarrow} \mathcal{O}_{\boldsymbol{E},s}^{2} 
           \stackrel{A}{\longrightarrow}\mathcal{O}_{\boldsymbol{E},s}^{2} 
           \longrightarrow \mathcal{O}_{\boldsymbol{E},s} 
           \longrightarrow \boldsymbol{k}(s) \longrightarrow 0
  $$
  where $AB=BA=f\cdot I_{2}$ is a reduced matrix factorisation of an equation
  $f$ of $\boldsymbol{E} \subset \mathbb{P}^{2}$. For example, if $s$ is a
  node, so that $\boldsymbol{E}$ is locally given by the polynomial $f = y^{2}
  - x^{3} -x^{2}\in\boldsymbol{k}[x,y]$, we can choose 
  $A=\bigl(\begin{smallmatrix}
    y&x^{2}+x\\x&y
  \end{smallmatrix}\bigr)$ and
  $B=\bigl(\begin{smallmatrix}
    y&-x^{2}-x\\-x&y
  \end{smallmatrix}\bigr)$ considered modulo $f$. More generally, any
  singular Weierstra{\ss} cubic $f$ can be written as $y\cdot y - R\cdot S$
  with $y, R,S$ all vanishing at the singular point. The off-diagonal elements
  of $A$ and $B$ are then formed by $\pm R,\pm S$. Therefore, all entries of
  the matrices $A$ and $B$ are elements of the maximal ideal of the local ring
  $\mathcal{O}_{\boldsymbol{E},s}$. Hence, the application of $\Hom(\ARG,
  \boldsymbol{k}(s))$ produces a complex with zero differential, which implies
  that $\Ext^{k}(\boldsymbol{k}(s), \boldsymbol{k}(s))$ is
  two-dimensional for all $k\ge 1$.
  In particular, $\Ext^{2}(\boldsymbol{k}(s), \boldsymbol{k}(s)) \cong
  \boldsymbol{k}^{2}$, and we can pick a non-zero element
  $w\in\Hom(\boldsymbol{k}(s), \boldsymbol{k}(s)[2])$. There exists a complex
  $X\in\Dbcoh(\boldsymbol{E})$ which sits in a distinguished triangle
  $$X\rightarrow \boldsymbol{k}(s) \stackrel{w}{\longrightarrow}
  \boldsymbol{k}(s)[2] \stackrel{+}{\longrightarrow}.$$
  Because the shift by one corresponds to multiplication by $-1$ in the
  Grothendieck group, this object $X$ is equal to zero in
  $\mathsf{K}(\boldsymbol{E})$. On the other hand, $X$ is
  indecomposable. Indeed, if $X$ would split, it must be $X\cong
  \boldsymbol{k}(s) \oplus \boldsymbol{k}(s)[1]$, because the only non-zero
  cohomology of $X$ is $H^{-1}(X) \cong \boldsymbol{k}(s)$ and $H^{0}(X) \cong
  \boldsymbol{k}(s)$. But, because $\Hom(\boldsymbol{k}(s)[1],
  \boldsymbol{k}(s)) = 0$, Lemma \ref{lem:PengXiao}, applied to the
  distinguished triangle
  \[\begin{CD}
    \boldsymbol{k}(s)[1] @>>> X @>>> \boldsymbol{k}(s) @>{+}>{w}>
  \end{CD}\]
  with $X\cong \boldsymbol{k}(s) \oplus \boldsymbol{k}(s)[1]$, implies $w=0$.
\end{example}

\begin{definition}[\cite{Stability}]
  A Harder-Narasimhan filtration (HNF) of an object $X \in
  \Dbcoh(\boldsymbol{E})$  is a finite collection of distinguished triangles
  \[
  \xymatrix@C=.5em
  {
    0\; \ar[rr] && F_{n}X \ar[rr] \ar[dl]_{\cong}&&   F_{n-1}X \ar[rr]
  \ar[dl]&& \dots \ar[rr] &&F_{1}X   \ar[rr]  &&  F_{0}X \ar@{=}[r]\ar[dl] & X
  \\
   & A_{n} \ar[lu]^{+} && A_{n-1} \ar[lu]^{+} &  & & & & &  A_0 \ar[lu]^{+}&
   \\
 }
 \]
 with $A_j\in\mathsf{P}(\varphi_j)$ and $A_{j}\ne 0$ for all $j$, such that
 $\varphi_{n} > \varphi_{n-1} > \cdots > \varphi_{0}.$
\end{definition}

If all ingredients of a HNF are shifted by one, we obtain a HNF of $X[1]$.
The shifted sheaves $A_{j}$ are called \emph{the semi-stable HN-factors} of
$X$ and we define $\varphi_{+}(X):=\varphi_{n}$ and
$\varphi_{-}(X):=\varphi_{0}$. Later, Theorem \ref{thm:uniqueHNF}, we show that
the HNF of an object $X$ is unique up to isomorphism. This justifies this
notation. For the moment, we keep in mind that $\varphi_{+}(X)$ and
$\varphi_{-}(X)$ might depend on the HNF and not only on the object $X$.

Before we proceed, we include a few remarks about the notation we use. 
Distinguished triangles in a triangulated category are either displayed in the
form
$X\rightarrow Y\rightarrow Z \stackrel{+}{\longrightarrow} 
  \quad\text{ or as }\quad
\xymatrix@C=.5em{
  X \ar[rr] && Y, \ar[dl]\\ & Z \ar[ul]^{+}
  }
$
where the arrow which is marked with $+$ is in fact a morphism $Z\rightarrow
X[1]$.

We shall use the octahedron axiom, the axiom (TR4) in Verdier's list, in the
following convenient form: if two morphisms $X\stackrel{f}{\longrightarrow} Y
\stackrel{g}{\longrightarrow} Z$  are given, for any three distinguished
triangles with bases $f, g$ and $g\circ f$ there exists a fourth distinguished
triangle which is indicated below by dashed arrows, such that we obtain the
following commutative diagram:

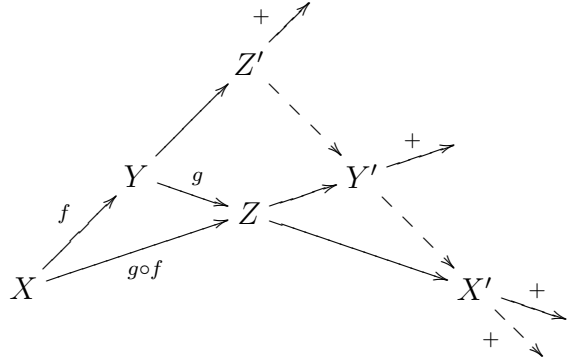
\begin{center}
  \mbox{\begin{xy} 0;<10mm,0mm>:0,  
      (0,3)      *+{Z'} ="top"      ,  
      (-3,0)     *+{X}  ="left"     ,  
      (3,0)      *+{X'} ="right"    ,  
      (-1.5,1.5) *+{Y}  ="midleft"  ,  
      (1.5,1.5)  *+{Y'} ="midright" ,  
      {"left";"midright":"right";"midleft",x} *+{Z}="center",
      {"left"  \ar@{->}^{f} "midleft"  \ar@{->}_{g\circ f} "center"},
      {"center" \ar@{->} "right" \ar@{->} "midright"},
      {"midleft"  \ar@{->}^{g} "center" \ar@{->} "top"},
      {"top" \ar@{-->} "midright"},
      {"midright" \ar@{-->} "right"},
      {"right" \ar@{-->}_{+} +(.9,-.9)},
      {"midright" \ar@{->}^{+} +"midright"-"center"},
      {"right" \ar@{->}^{+} +"center"-"midleft"},
      {"top" \ar@{->}^{+} +(.8,.8)}
  \end{xy}}  
\end{center}

The remainder of this section is devoted to the proofs of the crucial
properties of Harder-Narasimhan filtrations in triangulated categories. 
These properties can be found in \cite{Stability, GRK}, where most of
them appear to be either implicit or without a detailed proof.

\begin{lemma}\label{lem:connect}
  Let
  $\xymatrix@C=.4em{U \ar[rr]^{f} && X \ar[dl]\\ & V \ar[ul]^{+}}$
  and
  $A \longrightarrow V \longrightarrow V' \stackrel{+}{\longrightarrow}$
  be distinguished triangles. Then there exists a factorisation
  $U\longrightarrow W \stackrel{f'}{\longrightarrow} X$ 
  of $f$ and two distinguished triangles
  $$
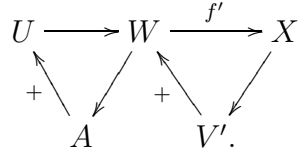
\xymatrix@C=.5em{U \ar[rr] && W \ar[dl]\ar[rr]^{f'} && X\ar[dl]\\
   & A \ar[ul]^{+} && V'.\ar[ul]^{+}}$$
\end{lemma}

\begin{proof}
  If we apply the octahedron axiom to the composition $A\rightarrow V
  \rightarrow U[1]$ we obtain the following commutative diagram, which gives
  the claim. 
  \begin{center}
  \mbox{\begin{xy} 0;<10mm,0mm>:0,  
      (0,3)      *+{V'}   ="top"      ,  
      (-3,0)     *+{A}    ="left"     ,  
      (3,0)      *+{X[1]} ="right"    ,  
      (-1.5,1.5) *+{V}    ="midleft"  ,  
      (1.5,1.5)  *+{W[1]} ="midright" ,  
      {"left";"midright":"right";"midleft",x} *+{U[1]}="center",
      {"left"  \ar@{->} "midleft"  \ar@{->} "center"},
      {"center" \ar@{->}_{f[1]} "right" \ar@{->} "midright"},
      {"midleft"  \ar@{->} "center" \ar@{->} "top"},
      {"top" \ar@{-->} "midright"},
      {"midright" \ar@{-->}^{f'[1]} "right"},
      {"right" \ar@{-->}_{+} +(.9,-.9)},
      {"midright" \ar@{->}^{+} +"midright"-"center"},
      {"right" \ar@{->}^{+} +"center"-"midleft"},
      {"top" \ar@{->}^{+} +(.8,.8)}
    \end{xy}}  
  \end{center}
\end{proof}

\begin{lemma}\label{lem:split}
  Let
  \[\xymatrix@C=.5em{
    0\; \ar[rr] && F_{n}V \ar[rr] \ar[dl]_{\cong}&&   F_{n-1}V \ar[rr]
  \ar[dl]&& \dots \ar[rr] &&F_{1}V   \ar[rr]  &&  F_{0}V \ar@{=}[r]\ar[dl] & V
  \\
   & A_{n} \ar[lu]^{+} && A_{n-1} \ar[lu]^{+} &  & & & & &  A_0 \ar[lu]^{+}&
  }\]
  be a HNF of $V\in\Dbcoh(\boldsymbol{E})$ and $F_{k}V \longrightarrow V
  \longrightarrow V' \stackrel{+}{\longrightarrow}$ a distinguished triangle
  with $1\le k \le n$. Then, $F_{k}V$ has a HNF with HN-factors $A_{n},
  A_{n-1}, \ldots, A_{k}$ and $V'$ one with HN-factors $A_{k-1},
  A_{k-2}, \ldots, A_{0}$. 
\end{lemma}

\begin{proof}
  The first statement is clear, because we can cut off the HNF of $V$ at
  $F_{k}V$ to obtain a HNF of $F_{k}V$. Let us define objects $F_{i}V'$ by
  exact triangles $F_{k}V \longrightarrow F_{i}V \longrightarrow F_{i}V'
  \stackrel{+}{\longrightarrow}$, where the first arrow is the composition of
  the morphisms in the HNF of $V$. Using the octahedron axiom, we obtain for
  any $i\le k$ a commutative diagram
  \begin{center}
  \mbox{\begin{xy} 0;<12mm,0mm>:0,  
      (0,3)      *+{F_{i}V'}   ="top"      ,  
      (-3,0)     *+{F_{k}V}    ="left"     ,  
      (3,0)      *+{A_{i-1}} ="right"    ,  
      (-1.5,1.5) *+{F_{i}V}    ="midleft"  ,  
      (1.5,1.5)  *+{F_{i-1}V'} ="midright" ,  
      {"left";"midright":"right";"midleft",x} *+{F_{i-1}V}="center",
      {"left"  \ar@{->} "midleft"  \ar@{->} "center"},
      {"center" \ar@{->} "right" \ar@{->} "midright"},
      {"midleft"  \ar@{->} "center" \ar@{->} "top"},
      {"top" \ar@{-->} "midright"},
      {"midright" \ar@{-->} "right"},
      {"right" \ar@{-->}_{+} +(.9,-.9)},
      {"midright" \ar@{->}^{+} +"midright"-"center"},
      {"right" \ar@{->}^{+} +"center"-"midleft"},
      {"top" \ar@{->}^{+} +(.8,.8)}
    \end{xy}}  
  \end{center}
  which implies the second claim.
\end{proof}

\begin{remark}\label{rem:split}
  The statement of Lemma \ref{lem:split} is true with identical proof if we
  relax the assumption of being a HNF by allowing $\varphi(A_{k}) =
  \varphi(A_{k-1})$ for the chosen value of $k$. 
\end{remark}

\begin{lemma}\label{lem:bounds}
  If
  \[\xymatrix@C=.5em{
    0\; \ar[rr] && F_{n}X \ar[rr] \ar[dl]_{\cong}&&   F_{n-1}X \ar[rr]
  \ar[dl]&& \dots \ar[rr] &&F_{1}X   \ar[rr]  &&  F_{0}X \ar@{=}[r]\ar[dl] & X
  \\
   & A_{n} \ar[lu]^{+} && A_{n-1} \ar[lu]^{+} &  & & & & &  A_0 \ar[lu]^{+}&
  }\]
  is a HNF of $X\in\Dbcoh(\boldsymbol{E})$ such that $A_{0}[k]$ is a sheaf,
  then $H^{k}(X)\ne 0$. In particular, the following implication is true:
  $$X\in \mathsf{D}^{\le m} \quad\Longrightarrow\quad
  \forall i\ge0: A_{i}\in  \mathsf{D}^{\le m}.$$
\end{lemma}

\begin{proof}
  The assumption $A_{0}[k]\in \Coh_{\boldsymbol{E}}$ means
  $H^{k}(A_{0})=A_{0}[k]\ne 0$ and $\varphi(A_{0}) \in (-k,-k+1]$. Because for
  all $i>0$ we have $\varphi(A_{i}) > \varphi(A_{0})$, we obtain
  $\varphi(A_{i})>-k$. This implies $H^{k+1}(A_{i})=0$ for all $i\ge 0$. The
  cohomology sequences of the distinguished triangles
  $F_{i+1}\longrightarrow F_{i}\longrightarrow A_{i}
  \stackrel{+}{\longrightarrow}$ imply $H^{k+1}(F_{i}X)=0$ for all $i>0$ and
  an exact sequence $H^{k}(X) \rightarrow H^{k}(A_{0}) \rightarrow
  H^{k+1}(F_{1}X)$, hence $H^{k}(X)\ne 0$. The statement about the other
  HN-factors $A_{i}$ follows now from $\varphi(A_{i})\ge \varphi(A_{0})$.
\end{proof}

\begin{proposition}
  Any non-zero object $X\in\Dbcoh(\boldsymbol{E})$ has a HNF.
\end{proposition}

\begin{proof}
  The existence of a  HNF for objects of $\Coh_{\boldsymbol{E}}$ is
  classically known, see \cite{HarderNarasimhan, Rudakov}. Therefore, we can
  proceed by induction on the number of non-zero cohomology sheaves of
  $X\in\Dbcoh(\boldsymbol{E})$. If $n$ is the largest integer with
  $H^{n}(X)\ne 0$, we have a distinguished triangle
  \begin{equation}
     \tau^{\le n-1} X \longrightarrow X \longrightarrow H^{n}(X)[-n]
    \stackrel{+}{\longrightarrow}
  \end{equation}
  By inductive hypothesis, there exists a HNF of $\tau^{\le n-1} X$. From
  Lemma \ref{lem:bounds} we conclude that all HN-factors of $\tau^{\le n-1} X$
  are in $\mathsf{D}^{\le n-1}$ and so $\varphi_{-}(\tau^{\le n-1} X)>-n+1$.

  Because $H^{n}(X)$ is a sheaf, we have $\varphi_{+}(H^{n}(X)[-n])
  \in (-n,-n+1]$, hence $\varphi_{-}(\tau^{\le n-1} X) >
  \varphi_{+}(H^{n}(X)[-n])$.

  We prove now for any distinguished triangle
  \begin{equation}
    \label{eq:induction}
    U \longrightarrow X \longrightarrow V \stackrel{+}{\longrightarrow}
  \end{equation}
  in which $V[n]$ is a coherent sheaf that the existence of a HNF for $U$ with
  $\varphi_{-}(U)> \varphi_{+}(V)$ implies the existence of a HNF of $X$.

  Because $V[n]$ is a sheaf, $V$ has a HNF and we proceed by induction on the
  number of HN-factors of $V$. Let $A$ be the leftmost object in a HNF of $V$,
  i.e.\/ $A\in\mathsf{P}(\varphi_{+}(V))$. By Lemma \ref{lem:connect} applied
  to the distinguished triangles (\ref{eq:induction}) and $A \longrightarrow V
  \longrightarrow V' \stackrel{+}{\longrightarrow}$, there exist two
  distinguished triangles in which $V'[n]$ is a coherent sheaf with a smaller
  number of HN-factors as $V$:
  $$\xymatrix@C=.5em{U \ar[rr] && W \ar[dl]\ar[rr] && X.\ar[dl]\\
    & A \ar[ul]^{+} && V'\ar[ul]^{+}}$$
  Because $\varphi_{-}(U)\ge \varphi(A) =\varphi_{+}(V)$, the left triangle can
  be concatenated to the given HNF of $U$ in order to provide a HNF for
  $W$. The start of the induction is covered as well: it is the case $V'=0$.
\end{proof}

\begin{lemma}\label{wesPT:ii}
  If $X,Y\in \Dbcoh(\boldsymbol{E})$ with $\varphi_{-}(X) > \varphi_{+}(Y)$,
  then $$\Hom(X,Y)=0.$$
\end{lemma}

\begin{proof}
  If $X,Y$ are semi-stable sheaves, this is well-known and follows easily from
  the definition of semi-stability. Because $\Hom(X,Y[k])=0$, if $X,Y$ are
  sheaves and $k<0$, the claim follows if $X\in \mathsf{P}(\varphi)$ and $Y\in
  \mathsf{P}(\psi)$ with $\varphi>\psi$. Let now $X\in\mathsf{P}(\varphi)$ and
  $Y\in \Dbcoh(\boldsymbol{E})$ with $\varphi > \varphi_{+}(Y)$. Let
  \[\xymatrix@C=.5em{
    0\; \ar[rr] && F_{m}Y \ar[rr] \ar[dl]_{\cong}&&   F_{m-1}Y \ar[rr]
  \ar[dl]&& \dots \ar[rr] &&F_{1}Y   \ar[rr]  &&  F_{0}Y \ar@{=}[r]\ar[dl] & Y
  \\
   & B_{m} \ar[lu]^{+} && B_{m-1} \ar[lu]^{+} &  & & & & &  B_0 \ar[lu]^{+}&
  }\]
  be a HNF of $Y$. We have $\varphi(B_{j})\le \varphi(B_{m}) =
  \varphi_{+}(Y)$, hence $\varphi(X)>\varphi(B_{j})$ and $\Hom(X,B_{j})=0$ for
  all $j$. If we apply the functor $\Hom(X,\ARG)$ to the distinguished 
  triangles $F_{j+1}Y \longrightarrow F_{j}Y \longrightarrow B_{j}
  \stackrel{+}{\longrightarrow}$, we obtain surjections $\Hom(X,F_{j+1}Y)
  \twoheadrightarrow \Hom(X,F_{j}Y)$. From $\Hom(X,F_{m}Y)=
  \Hom(X,B_{m})=0$, we obtain $\Hom(X,Y)=\Hom(X,F_{0}Y)=0$.

  Let now $X,Y$ be arbitrary non-zero objects of $\Dbcoh(\boldsymbol{E})$
  which satisfy $\varphi_{-}(X) > \varphi_{+}(Y)$. If
  \[\xymatrix@C=.5em{
    0\; \ar[rr] && F_{n}X \ar[rr] \ar[dl]_{\cong}&&   F_{n-1}X \ar[rr]
  \ar[dl]&& \dots \ar[rr] &&F_{1}X   \ar[rr]  &&  F_{0}X \ar@{=}[r]\ar[dl] & X
  \\
   & A_{n} \ar[lu]^{+} && A_{n-1} \ar[lu]^{+} &  & & & & &  A_0 \ar[lu]^{+}&
 }\]
 is a HNF of $X$, we have $\varphi(A_{i})\ge \varphi(A_{0})=\varphi_{-}(X) >
  \varphi_{+}(Y)$. We know already $\Hom(A_{i},Y)=0$ for all $i\ge 0$. If we
  apply the functor $\Hom(\ARG,Y)$ to the distinguished triangles
  $F_{i+1}X \longrightarrow F_{i}X \longrightarrow A_{i}
  \stackrel{+}{\longrightarrow}$, we obtain injections $\Hom(F_{i}X,Y)
  \hookrightarrow \Hom(F_{i+1}X,Y)$. Again, this implies $\Hom(X,Y)=0$.
\end{proof}

\begin{theorem}[\cite{Stability,GRK}]\label{thm:uniqueHNF}
  The HNF of any non-zero object $X\in\Dbcoh(\boldsymbol{E})$ is unique up to
  unique isomorphism.
\end{theorem}

\begin{proof}
  If
  \[\xymatrix@C=.5em{
    0\; \ar[rr] && F_{n}X \ar[rr] \ar[dl]_{\cong}&&   F_{n-1}X \ar[rr]
    \ar[dl]&& \dots \ar[rr] &&F_{1}X   \ar[rr]  &&  F_{0}X \ar@{=}[r]\ar[dl] &
    X\\
    & A_{n} \ar[lu]^{+} && A_{n-1} \ar[lu]^{+} &  & & & & &  A_0 \ar[lu]^{+}&
  }\]
  and
  \[\xymatrix@C=.5em{
    0\; \ar[rr] && G_{m}X \ar[rr] \ar[dl]_{\cong}&&   G_{m-1}X \ar[rr]
    \ar[dl]&& \dots \ar[rr] &&G_{1}X   \ar[rr]  &&  G_{0}X \ar@{=}[r]\ar[dl] &
    X\\
    & B_{m} \ar[lu]^{+} && B_{m-1} \ar[lu]^{+} &  & & & & &  B_0 \ar[lu]^{+}&
  }\]
  are HNFs of $X$, we have to show that there exist unique isomorphisms of
  distinguished triangles for any $k\ge 0$
  \[\begin{CD}
    F_{k+1}X      @>>> F_{k}X      @>>> A_{k}       @>{+}>>\\
    @VV{f_{k+1}}V      @VV{f_{k}}V      @VV{g_{k}}V \\
    G_{k+1}X      @>>> G_{k}X      @>>> B_{k}       @>{+}>>
  \end{CD}\]
  with $f_{0}=\mathsf{Id}_{X}$. This is obtained by induction on $k\ge0$ from
  the following claim: if an isomorphism $f:F\rightarrow G$ and two
  distinguished triangles
  $F' \longrightarrow F \longrightarrow A \stackrel{+}{\longrightarrow }$ and
  $G' \longrightarrow G \longrightarrow B \stackrel{+}{\longrightarrow }$ are
  given such that $A\in\mathsf{P}(\varphi), B\in\mathsf{P}(\psi)$ and $F',G'$
  have HNFs with $\varphi_{-}(F')>\varphi$ and $\varphi_{-}(G')>\psi$, then
  there exist unique isomorphisms $f':F'\rightarrow G'$ and $g:A\rightarrow B$
  such that $(f',f,g)$ is a morphism of triangles. In particular,
  $\varphi=\psi$.

  Without loss of generality, we may assume $\varphi\ge \psi$. This implies
  $\varphi_{-}(F'[1]) > \varphi_{-}(F') > \psi$. Lemma \ref{wesPT:ii} implies
  therefore $\Hom(F',B) = \Hom(F'[1],B) = 0$. From \cite{Asterisque100},
  Proposition 1.1.9, we obtain the existence and uniqueness of the morphisms
  $f',g$. It remains to show that they are isomorphisms. If $g$ were zero, the
  second morphism in the triangle $G' \longrightarrow G
  \stackrel{0}{\longrightarrow} B \stackrel{+}{\longrightarrow }$ would be
  zero. Hence, $B$ were a direct summand of  $G'[1]$ which implies
  $\Hom(G'[1],B)\ne 0$. This contradicts Lemma \ref{wesPT:ii}, because
  $\varphi_{-}(G'[1]) > \varphi(G') > \psi=\varphi(B)$. Hence, $g\ne 0$ and
  Lemma \ref{wesPT:ii} implies $\varphi(A)\le \varphi(B)$, i.e.\/
  $\varphi=\psi$. So, the same reasoning as before gives a unique morphism of
  distinguished triangles in the other direction. The composition of both are
  the respective identities of $F' \longrightarrow F \longrightarrow A
  \stackrel{+}{\longrightarrow }$ and $G' \longrightarrow G \longrightarrow B
  \stackrel{+}{\longrightarrow }$ respectively, which follows again from the
  uniqueness part of \cite{Asterisque100}, Proposition 1.1.9. This proves the
  claim. 
\end{proof}

We need the following useful lemma.

\begin{lemma}(\cite{PengXiao}, Lemma 2.5)\label{lem:PengXiao}
Let $\mathsf{D}$ be a triangulated category and 
\[\begin{CD}
 F  @>>> G  @>>> H_1 \oplus H_2 @>{+}>{(0,w)}>
\end{CD}\]
be a distinguished triangle in $\mathsf{D}$. Then $G\cong H_{1}\oplus G'$
splits and the given triangle is isomorphic to  
\[\begin{CD}
F  
@>{\bigl(\begin{smallmatrix} 0\\g \end{smallmatrix}\bigr)}>> 
H_1 \oplus G'   
@>{\bigl(\begin{smallmatrix}1&0\\0&f' \end{smallmatrix}\bigr)}>>
H_1 \oplus H_2 
@>{+}>{(0,w)}>
\end{CD}\]
Dually, if 
\[\begin{CD}
 F  @>{\bigl(\begin{smallmatrix} 0\\g \end{smallmatrix}\bigr)}>> 
 G_{1}\oplus G_{2}  @>>> H @>{+}>>
\end{CD}\]
is a distinguished triangle then $H\cong G_{1}\oplus H'$ and the given
triangle is isomorphic to
\[\begin{CD}
F  
@>{\bigl(\begin{smallmatrix} 0\\g \end{smallmatrix}\bigr)}>> 
G_1 \oplus G_{2}
@>{\bigl(\begin{smallmatrix}1&0\\0&f' \end{smallmatrix}\bigr)}>>
G_1 \oplus H'
@>{+}>{(0,w)}>
\end{CD}\]
\end{lemma} 

The results in this section are true for more general triangulated
categories than $\Dbcoh(\boldsymbol{E})$. Without changes, the proofs apply if
we replace $\Dbcoh(\boldsymbol{E})$ by the bounded derived category of an
Abelian category which is equipped with the notion of stability in the sense
of \cite{Rudakov}. In particular, these results hold for polynomial stability
on the triangulated categories $\Dbcoh(X)$ where $X$ is a projective 
variety over $\boldsymbol{k}$.

\section{The structure of the bounded derived category of coherent sheaves on
  a singular Weiersta{\ss} curve}\label{sec:dercat}

In this section, we prove the main results on which our understanding of
$\Dbcoh(\boldsymbol{E})$ is based. Again, $\boldsymbol{E}$ denotes a
Weierstra{\ss} curve. Our main focus is on the singular case, however all the
results remain true in the smooth case as well. 
A speciality of this category is the
non-vanishing result Proposition \ref{wesPT}. Unlike the smooth case, there
exist indecomposable objects in $\Dbcoh(\boldsymbol{E})$, which are not
semi-stable. Their Harder-Narasimhan factors are characterised in Proposition
\ref{prop:extreme}. We propose to visualise indecomposable objects by their
``shadows''. As an application of our results, we give a complete
characterisation of all spherical objects in $\Dbcoh(\boldsymbol{E})$. As a
consequence, we show that the group of exact auto-equivalences acts
transitively on the set of spherical objects. This answers a question which was
posed by Polishchuk \cite{YangBaxter}.

Let us set up some notation. 
For any $\varphi\in(0,1]$ we denote by $\mathsf{P}(\varphi)^{s} \subset
\mathsf{P}(\varphi)$ the full subcategory of stable sheaves with phase
$\varphi$. We extend this definition to all $\varphi\in\mathbb{R}$ by
requiring $\mathsf{P}(\varphi +n)^{s} = \mathsf{P}(\varphi)^{s}[n]$ for all
$n\in\mathbb{Z}$ and all $\varphi\in\mathbb{R}$.

We already know the structure of $\mathsf{P}(1)^{s}$. Because $\mathsf{P}(1)$
is the category of coherent torsion sheaves on $\boldsymbol{E}$, the objects
of $\mathsf{P}(1)^{s}$ are precisely the structure sheaves $\boldsymbol{k}(x)$
of closed points $x\in\boldsymbol{E}$. In order to understand the structure of
all the other categories $\mathsf{P}(\varphi)^{s}$, we use Fourier-Mukai
transforms. Our main technical tool will be the transform $\mathbb{F}$ which
was studied in \cite{BurbanKreussler}. It depends on the choice of a regular
point $p_{0}\in\boldsymbol{E}$.
Let us briefly recall its definition and main properties. It was defined
with the aid of Seidel-Thomas twists \cite{SeidelThomas}, which are functors
$T_{E}: \Dbcoh(\boldsymbol{E}) \rightarrow \Dbcoh(\boldsymbol{E})$
depending on a spherical object $E\in\Dbcoh(\boldsymbol{E})$. On objects, these
functors are characterised by the existence of a distinguished triangle
$$\boldsymbol{R}\Hom(E,F) \otimes E \rightarrow F \rightarrow T_{E}(F)
\stackrel{+}{\longrightarrow}.$$ 
If $p_{0}\in\boldsymbol{E}$ is a smooth point, the functor
$T_{\boldsymbol{k}(p_{0})}$ is isomorphic to the tensor product with the
locally free sheaf $\mathcal{O}_{\boldsymbol{E}}(p_{0})$, see
\cite{SeidelThomas}, 3.11. We defined 
$$\mathbb{F} := 
T_{\boldsymbol{k}(p_{0})}T_{\mathcal{O}}T_{\boldsymbol{k}(p_{0})}.$$
In \cite{SeidelThomas} is was shown that twist functors can be described as
integral transforms and that $\mathbb{F}$ is isomorphic to the functor
$\FM^{\mathcal{P}}$, which is given by 
$$\FM^{\mathcal{P}}(\ARG) := 
\boldsymbol{R}\pi_{2\ast}(\mathcal{P}\dtens \pi_{1}^{\ast}(\ARG)),$$
where $\mathcal{P}=\mathcal{I}_{\Delta}\otimes
\pi_{1}^{\ast}\mathcal{O}(p_{0}) \otimes
\pi_{2}^{\ast}\mathcal{O}(p_{0})[1]$. This is a shift of a coherent sheaf on
$\boldsymbol{E}\times \boldsymbol{E}$, on which we denote the ideal of the
diagonal by $\mathcal{I}_{\Delta} \subset
\mathcal{O}_{\boldsymbol{E}\times\boldsymbol{E}}$ and the two projections by
$\pi_{1}, \pi_{2}$. 

In order to understand the effect of $\mathbb{F}$ on rank and degree, we look
at the distinguished triangle 
$$\boldsymbol{R}\Hom(\mathcal{O},F) \otimes \mathcal{O} \rightarrow F
\rightarrow T_{\mathcal{O}}(F) \stackrel{+}{\longrightarrow}.$$ 
The additivity of rank and degree implies $\rk(T_{\mathcal{O}}(F))= \rk(F) -
\deg(F)$ and 
$\deg(T_{\mathcal{O}}(F))= \deg(F)$. On the other hand, it is well-known that 
$\deg(T_{\boldsymbol{k}(p_{0})}(F)) = \deg(F)+\rk(F)$ and 
$\rk(T_{\boldsymbol{k}(p_{0})}(F)) = \rk(F)$.
So, if we use $[\mathcal{O}_{\boldsymbol{E}}], -[\boldsymbol{k}(p_{0})]$ as a
basis of $\mathsf{K}(\boldsymbol{E})$, which means that we use coordinates
$(\rk,-\deg)$, then the action of $T_{\mathcal{O}}, T_{\boldsymbol{k}(p_{0})}$
and of $\mathbb{F}$ on $\mathsf{K}(\boldsymbol{E})$ is given by the matrices
  $$ 
  \begin{pmatrix}
    1&1\\0&1
  \end{pmatrix}, 
  \begin{pmatrix}
    1&0\\-1&1
  \end{pmatrix} \quad \;\text{and}\quad
  \begin{pmatrix}
    0&1\\-1&0
  \end{pmatrix}\;\text{respectively.}
  $$
In particular, for any object $F\in\Dbcoh(\boldsymbol{E})$ which has a slope,
we have $\mu(T_{\boldsymbol{k}(p_{0})}(F)) = \mu(F)+1$ and 
$\mu(\mathbb{F}(F))=-\frac{1}{\mu(F)}$ using the usual conventions in dealing
with $\infty$. 

If $F$ is a sheaf or a twist thereof, we defined the phase
$\varphi(F)$. In order to understand the effect of $\mathbb{F}$ on phases, it
is not sufficient to know its effect on the slope. This is because the slope
determines the phase modulo $2\mathbb{Z}$ only.
However, if $F$ is a coherent sheaf, the description of $\mathbb{F}$ as
$\FM^{\mathcal{P}}$ shows that $\mathbb{F}(F)$ can have non-vanishing
cohomology in degrees $-1$ and $0$ only. If, in addition, $\mathbb{F}(F)$ is a
shifted sheaf, this implies $\varphi(\mathbb{F}(F))\in (0,2]$.
From the formula for the slope it is now clear that $\varphi(\mathbb{F}(F)) =
\varphi(F)+\frac{1}{2}$ for any shifted coherent sheaf $F$.

The following result was first shown in \cite{Nachr}. We give an independent
proof here, which was inspired by \cite{Nachr}, Lemma 3.1. 

\begin{theorem}\label{thm:mother}
  $\mathbb{F}$ sends semi-stable sheaves to semi-stable sheaves.
\end{theorem}

\begin{proof}
  Note that, by definition, a semi-stable sheaf of positive rank is
  automatically torsion free. The only sheaf with degree and rank equal to
  zero is the zero sheaf. Throughout this proof, we let $\mathcal{F}$ be a
  semi-stable sheaf on $\boldsymbol{E}$.
  If $\deg(\mathcal{F})=0$ this sheaf is torsion free and the claim
  was shown in \cite{BurbanKreussler}, Thm.\/ 2.21, see also \cite{FMmin}.
  For the sake of clarity we would like to stress here the fact that
  \cite{BurbanKreussler}, Section 2, deals with nodal as well as cuspidal
  Weierstra{\ss} curves. 

  Next, suppose $\deg(\mathcal{F})>0$. If
  $\rk(\mathcal{F})=0$, $\mathcal{F}$ is a coherent torsion sheaf. Again, the
  claim follows from \cite{BurbanKreussler}, Thm.\/ 2.21 and Thm.\/ 2.18,
  where it was shown that $\mathbb{F}\circ\mathbb{F}= i^{\ast}[1]$, for any
  Weierstra{\ss} curve. Here, $i:\boldsymbol{E} \rightarrow \boldsymbol{E}$ is
  the involution which fixes the singularity and which corresponds to taking
  the inverse on the smooth part of $\boldsymbol{E}$ with its group structure
  in which $p_{0}$ is the neutral element.

  Therefore, we may suppose $\mathcal{F}$ is torsion free. As observed before,
  the complex $\mathbb{F}(\mathcal{F})\in\Dbcoh(\boldsymbol{E})$ can have
  non-vanishing cohomology in degrees $-1$ and $0$ only. We are going to show
  that $\mathbb{F}(\mathcal{F})[-1]$ is a sheaf, which is equivalent to the
  vanishing of the cohomology object
  $\mathcal{H}^{0}(\mathbb{F}(\mathcal{F}))\in\Coh_{\boldsymbol{E}}$.
  Recall from \cite {BurbanKreussler}, Lemma 2.13, that for any smooth point
  $x\in\boldsymbol{E}$ the sheaf of degree zero $\mathcal{O}(x-p_{0})$
  satisfies  
  $\mathbb{F}(\mathcal{O}(x-p_{0})) \cong T_{\mathcal{O}}(\mathcal{O}(x))
  \cong \boldsymbol{k}(x)$. Moreover, if $s\in\boldsymbol{E}$ denotes the
  singular point, $n:\mathbb{P}^{1}\rightarrow \boldsymbol{E}$ 
  the normalisation and
  $\widetilde{\mathcal{O}}:=n_{\ast}(\mathcal{O}_{\mathbb{P}^{1}})$, then
  $\mathbb{F}(\widetilde{\mathcal{O}}(-p_{0})) \cong
  T_{\mathcal{O}}(\widetilde{\mathcal{O}}) \cong \boldsymbol{k}(s)$. The sheaf
  $\widetilde{\mathcal{O}}(-p_{0})$ has degree zero on $\boldsymbol{E}$.
  Because $\mathbb{F}$ is an equivalence, we obtain isomorphisms
  \begin{align*}
    \Hom(\mathbb{F}(\mathcal{F}),\boldsymbol{k}(x)) &\cong
    \Hom(\mathcal{F}, \mathcal{O}(x-p_{0}))\\
    \intertext{and}
    \Hom(\mathbb{F}(\mathcal{F}),\boldsymbol{k}(s)) &\cong
    \Hom(\mathcal{F}, \widetilde{\mathcal{O}}(-p_{0}))
  \end{align*}
  where $x\in\boldsymbol{E}$ is an arbitrary smooth point. These vector
  spaces vanish as $\mathcal{F}$ was assumed to be semi-stable and of positive
  degree.
  
  Because cohomology of the complex $\mathbb{F}(\mathcal{F})$ vanishes in
  positive degree, there is a canonical 
  morphism $\mathbb{F}(\mathcal{F})\rightarrow
  \mathcal{H}^{0}(\mathbb{F}(\mathcal{F}))$ in $\Dbcoh(\boldsymbol{E})$, which
  induces an injection of functors
  $\Hom(\mathcal{H}^{0}(\mathbb{F}(\mathcal{F})), \ARG) \hookrightarrow 
  \Hom(\mathbb{F}(\mathcal{F}), \ARG)$. Therefore, the vanishing which was
  obtained above, shows
  $$\Hom(\mathcal{H}^{0}(\mathbb{F}(\mathcal{F})), \boldsymbol{k}(y)) = 0$$
  for any point $y\in\boldsymbol{E}$. This implies the vanishing of the sheaf
  $\mathcal{H}^{0}(\mathbb{F}(\mathcal{F}))$. Hence, 
  $\widehat{\mathcal{F}}:=\mathbb{F}(\mathcal{F})[-1]$ is a coherent sheaf and
  the definition of $\mathbb{F}$ implies that there is an exact sequence of
  coherent sheaves
  $$0\rightarrow \widehat{\mathcal{F}}(-p_{0}) \rightarrow
  H^{0}(\mathcal{F}(p_{0})) \otimes \mathcal{O}_{\boldsymbol{E}} \rightarrow
  \mathcal{F}(p_{0}) \rightarrow 0.$$ 
  This sequence implies, in particular, that $\widehat{\mathcal{F}}$ is
  torsion free. 

  Before we proceed to show that $\widehat{\mathcal{F}}$ is semi-stable, we
  apply duality to prove that $\mathbb{F}(\mathcal{F})$ is a sheaf if
  $\deg(\mathcal{F})<0$. Let us denote the dualising functor by $\mathbb{D}:=
  \boldsymbol{R}\mathcal{H}om(\ARG, \mathcal{O}_{\boldsymbol{E}})$. This
  functor satisfies $\mathbb{D}\mathbb{D}\cong \boldsymbol{1}$. In
  \cite{BurbanKreusslerRel}, Cor.\/ 3.4, we have shown that there exists an
  isomorphism
  $$\mathbb{D}\mathbb{F} [-1] \cong i^{\ast} \mathbb{F} \mathbb{D}.$$
  Using $\mathbb{D}\circ[1] \cong [-1]\circ \mathbb{D}$, this implies 
  $$\mathbb{F} \cong \mathbb{D}i^{\ast}[-1]\mathbb{F}\mathbb{D}.$$
  Because $\mathcal{F}$ is a torsion free sheaf on a curve, it is
  Cohen-Macaulay and since $\boldsymbol{E}$ is Gorenstein, this implies
  $\mathcal{E}xt^{i}(\mathcal{F},\mathcal{O}) = 0$ for any $i>0$.  
  Therefore, we have $\mathbb{D}(\mathcal{F})\cong \mathcal{F}^{\vee}$ and
  this is a semi-stable coherent sheaf of positive degree. Thus, $[-1]\circ
  \mathbb{F}$ sends $\mathcal{F}^{\vee}$ to a torsion free sheaf, on which
  $\mathbb{D}$ is just the usual dual. Now, we see that
  $\mathbb{F}(\mathcal{F})$ is a torsion free sheaf if $\mathcal{F}$ was 
  semi-stable and of negative degree. 

  It remains to prove that $\mathbb{F}$ preserves semi-stability. If
  $\deg(\mathcal{F})=0$ or $\mathcal{F}$ is a torsion sheaf, this was shown
  for any Weierstra{\ss} curve in \cite{BurbanKreussler}.
  If $\deg(\mathcal{F})\ne 0$ the proof is based upon
  $\mathbb{F}\mathbb{F}[-1]\cong i^{\ast}$, see \cite{BurbanKreussler}, Thm.\/
  2.18. Suppose $\deg(\mathcal{F})>0$, then
  $\mathbb{F}(\widehat{\mathcal{F}})\cong i^{\ast}(\mathcal{F})$ and this is a
  coherent sheaf. 
  If $\widehat{\mathcal{F}}$ were not semi-stable, there would exist a
  semi-stable sheaf $\mathcal{G}$ with $\mu(\widehat{\mathcal{F}}) >
  \mu(\mathcal{G})$ and a non-zero morphism $\widehat{\mathcal{F}} \rightarrow
  \mathcal{G}$. Because $\mu(\widehat{\mathcal{F}}) = -1/\mu(\mathcal{F})<0$, 
  $\mathbb{F}(\mathcal{G})$ is a coherent sheaf and application of
  $\mathbb{F}$ produces a non-zero morphism $i^{\ast}(\mathcal{F}) \cong
  \mathbb{F}(\widehat{\mathcal{F}}) \rightarrow
  \mathbb{F}(\mathcal{G})$. However, $\mu(i^{\ast}(\mathcal{F})) =
  \mu(\mathcal{F}) > -1/\mu(\mathcal{G}) =  \mu(\mathbb{F}(\mathcal{G}))$
  contradicts semi-stability of $i^{\ast}(\mathcal{F})$. Hence,
  $\widehat{\mathcal{F}}$ is semi-stable. The proof in the case
  $\deg(\mathcal{F})<0$ starts with a non-zero morphism
  $\mathcal{U}\rightarrow \mathbb{F}(\mathcal{F})$ and proceeds similarly.  
\end{proof}

It was shown in \cite{BurbanKreussler} that we obtain an action of the group
$\widetilde{\SL}(2,\mathbb{Z})$ on $\Dbcoh(\boldsymbol{E})$ by sending
generators of this group to $T_{\mathcal{O}}$, $T_{\boldsymbol{k}(p_{0})}$ and
the translation functor $[1]$ respectively.
Let us denote $$\mathsf{Q}:=\{\varphi\in\mathbb{R}\mid
\mathsf{P}(\varphi) \text{ contains a non-zero object}\}.$$
The action of a group $G$ on $\mathsf{Q}$ is called \emph{monotone}, if
$\varphi\le\psi$ implies $g\cdot\varphi\le g\cdot\psi$ for every $g\in G$ and
$\varphi,\psi\in \mathsf{Q}$.

\begin{proposition}\label{prop:transit}
  The $\widetilde{\SL}(2,\mathbb{Z})$-action on $\Dbcoh(\boldsymbol{E})$
  induces a monotone and transitive action on the set $\mathsf{Q}$. All
  isotropy groups of this action are isomorphic to $\mathbb{Z}$.
\end{proposition}

\begin{proof}
  As seen above, for any $\psi\in\mathsf{Q}$ and $0\ne
  A\in\mathsf{P}(\psi)$, we have $\varphi(\mathbb{F}(A)) =
  \varphi(A)+\frac{1}{2}$ and $\mu(T_{\boldsymbol{k}(p_{0})}(A)) = \mu(A)+1$. 
  Therefore, by Theorem \ref{thm:mother} it is clear that we obtain an induced
  monotone action of $\widetilde{\SL}(2,\mathbb{Z})$ on $\mathsf{Q}$.
  The group $\SL(2,\mathbb{Z})$ acts transitively on the set of all pairs of
  co-prime integers which we interpret as primitive vectors of the lattice
  $\mathbb{Z}\oplus i\mathbb{Z}\subset\mathbb{C}$. Hence, the action 
  of $\widetilde{\SL}(2,\mathbb{Z})$ on $\mathsf{Q}$ is transitive as
  well. So, all isotropy groups are isomorphic. Finally, it is easy
  to see that the isotropy group of $1\in\mathsf{Q}$ is generated by
  $T_{\boldsymbol{k}(p_{0})}$.
\end{proof}

As an important consequence we obtain the following clear structure result
for the slices $\mathsf{P}(\varphi)$. 

\begin{corollary}\label{cor:equiv}
  The category $\mathsf{P}(\varphi)$ of semi-stable objects of phase
  $\varphi\in\mathsf{Q}$ is equivalent to the  category $\mathsf{P}(1)$ of
  torsion sheaves. Any such equivalence restricts to an equivalence between
  $\mathsf{P}(\varphi)^{s}$ and $\mathsf{P}(1)^{s}$. Under such an
  equivalence, stable vector bundles correspond to  structure sheaves of
  smooth points. Moreover, if $\varphi\in(0,1)\cap \mathsf{Q}$,
  $\mathsf{P}(\varphi)^{s}$ contains a unique torsion free sheaf, which is not
  locally free. It correspond to the structure sheaf
  $\boldsymbol{k}(s)\in\mathsf{P}(1)^{s}$ of the singular point.  
\end{corollary}

Recall that an object $E\in\Dbcoh(\boldsymbol{E})$ is called \emph{perfect},
if it is isomorphic in the derived category to a bounded complex of locally
free sheaves of finite rank. Thus, a sheaf or shift thereof is called perfect,
if it is perfect as an object in $\Dbcoh(\boldsymbol{E})$. 
If $\boldsymbol{E}$ is smooth, any object in $\Dbcoh(\boldsymbol{E})$ is
perfect. However, if $s\in\boldsymbol{E}$ is a singular point, the torsion
sheaf $\boldsymbol{k}(s)$ is not perfect. 

If $\boldsymbol{E}$ is singular with one singularity $s\in\boldsymbol{E}$, the
category $\mathsf{P}(1)^{s}$ contains precisely one object which is not
perfect, the object $\boldsymbol{k}(s)$. 
Hence, by Proposition \ref{prop:transit}, for any $\varphi\in\mathsf{Q}$ there
is precisely one element in $\mathsf{P}(\varphi)^{s}$ which is not perfect. We
shall refer to it as the \emph{extreme} stable element with phase
$\varphi$. So, the sheaf $\boldsymbol{k}(s)$ is the extreme stable element
with phase $1$. The extreme stable element is never locally free. A stable
object is either perfect or extreme.

We shall need the following version of Serre duality, which can be deduced
easily from standard versions:

If $E,F\in\Dbcoh(\boldsymbol{E})$ and at least one of them is perfect, then
there is a bi-functorial isomorphism
\begin{equation}
  \label{wesPT:i}\Hom(E,F) \cong \Hom(F,E[1])^{\ast}.
\end{equation}

If neither of the objects is perfect, this is no longer true. For example,
$\Hom(\boldsymbol{k}(s),\boldsymbol{k}(s))\cong \boldsymbol{k}$, but
$\Hom(\boldsymbol{k}(s),\boldsymbol{k}(s)[1]) \cong
\Ext^{1}(\boldsymbol{k}(s),\boldsymbol{k}(s)) \cong \boldsymbol{k}^{2}$.

Any object $X$ in the Abelian category $\mathsf{P}(\varphi)$ has a
Jordan-H\"older filtration (JHF) 
$$0\subset F_{n}X \subset \ldots \subset F_{1}X \subset F_{0}X = X$$
with stable JH-factors $J_{i}=F_{i}X/F_{i+1}X \in
\mathsf{P}(\varphi)^{s}$. The graded object $\oplus_{i=0}^{n}J_{i}$ is
determined by $X$. Observe that for any two objects $J\not\cong
J'\in\mathsf{P}(\varphi)^{s}$ we can apply Serre duality because at most one
of them is non-perfect.

\begin{corollary}\label{cor:sheaves}
  \begin{enumerate}
  \item\label{cor:i} If $\varphi,\psi \in \mathsf{Q}$ with $\varphi -1 < \psi
  \le \varphi$ there exists $\Phi\in \widetilde{\SL}(2,\mathbb{Z})$ such that
  $\Phi(\varphi)=1$ and $\Phi(\psi)\in(0,1]$.
  \item\label{cor:ii} If $A,B\in\mathsf{P}(\varphi)^{s}$, then $A\cong B \iff
  \Hom(A,B)\ne 0.$ 
  \item\label{cor:iii} If $0\ne X\in\mathsf{P}(\varphi)$ and $0\ne
  Y\in\mathsf{P}(\psi)$ with $\varphi < \psi < \varphi+1$, then $\Hom(X,Y)\ne
  0$. 
  \item \label{cor:iv} If $J\in\mathsf{P}(\varphi)^{s}$ is not a JH-factor of 
  $X\in \mathsf{P}(\varphi)$, for all $i\in\mathbb{Z}$ we have
  $\Hom(J,X[i])=0$.
  \item \label{cor:v} If $X\in\mathsf{P}(\varphi)$ is indecomposable, all its
    JH-factors are isomorphic to each other. 
  \item \label{cor:vi} If $X,Y\in\mathsf{P}(\varphi)$ are non-zero
    indecomposable objects, both with the same JH-factor, then $\Hom(X,Y) \ne
    0$. 
  \end{enumerate}
\end{corollary}

\begin{proof}
  (\ref{cor:i}) This follows from Proposition \ref{prop:transit} because the
  shift functor corresponds to an element in the centre of
  $\widetilde{\SL}(2,\mathbb{Z})$ and therefore $\Phi(\mathsf{P}(\varphi)) =
  \mathsf{P}(1)$ implies $\Phi(\mathsf{P}(\varphi-1)) = \mathsf{P}(0)$.

  (\ref{cor:ii}) The statement is clear in case $\varphi=1$ and follows from
  (\ref{cor:i}) in the general case.

  (\ref{cor:iii}) Using (\ref{cor:i}) we can assume $\psi=1$, which means
  that $Y$ is a coherent torsion sheaf. By Proposition \ref{prop:transit} this
  implies $\varphi\in(0,1)$ and $X$ is a torsion free coherent sheaf. If
  $Y\in\mathsf{P}(1)^{s}$ the statement is clear, because any torsion free
  sheaf has a non-zero morphism to any $Y=\boldsymbol{k}(x)$,
  $x\in\boldsymbol{E}$. If $Y\in\mathsf{P}(1)$ is arbitrary, there exists a
  point $x\in\boldsymbol{E}$ and a non-zero morphism $\boldsymbol{k}(x)
  \rightarrow Y$. The claim follows now from left-exactness of the functor
  $\Hom(X,\ARG)$.

  (\ref{cor:iv}) If $J'\in\mathsf{P}(\varphi)^{s}$ is a JH-factor of $X$, we
  have $J\not\cong J'$. From (\ref{cor:ii}) and Serre duality together with
  Lemma \ref{wesPT:ii} we obtain $\Hom(J,J'[i])=0$ for any
  $i\in\mathbb{Z}$. Using the JHF of $X$, the claim now follows.

  (\ref{cor:v}) It is easy to prove by induction that any
  $X\in\mathsf{P}(\varphi)$ can be split as a finite direct sum $X\cong\oplus
  X_{k}$, where each $X_{k}$ has all JH-factors isomorphic to a single element
  $J_{k}\in\mathsf{P}(\varphi)^{s}$. This implies, the claim.

  (\ref{cor:vi}) By (\ref{cor:i}) we may assume
  $\varphi(X)=\varphi(Y)=1$. This means, both objects are indecomposable
  torsion sheaves with support at the singular point $s\in\boldsymbol{E}$.
  Such sheaves always have an epimorphism to and a monomorphism from the
  extreme object $\boldsymbol{k}(s)$, hence the claim.
\end{proof}

It is interesting and important to note that an indecomposable semi-stable
object can be perfect even though all its JH-factors are extreme. This is
made explicit in \cite{BurbanKreussler}, Section 4, in the case of the
category $\mathsf{P}(1)$ of coherent torsion sheaves. If $\boldsymbol{E}$ is
nodal, there are two kinds of indecomposable torsion sheaves with support at
the node $s\in\boldsymbol{E}$: the so-called \emph{bands} and
\emph{strings}. The bands are perfect, whereas the strings are not
perfect. Using the action of $\widetilde{\SL}(2,\mathbb{Z})$ this carries over
to all other categories $\mathsf{P}(\varphi)$ with $\varphi\in\mathsf{Q}$. 

An object $X\in\mathsf{P}(\varphi)$ will be called \emph{extreme} if it
does not have a direct summand which is perfect. This implies that, but is not
equivalent to the property that all its JH-factors are extreme. An example can
be found below, see Ex.~\ref{ex:extremefactors}.
From the above we deduce that any $X\in\mathsf{P}(\varphi)$ can be split as a
direct sum $X\cong X^{e}\oplus X^{p}$ with $X^{e}$ extreme and $X^{p}$
perfect. All direct summands of the extreme part have the unique
extreme stable element with phase $\varphi$ as its JH-factors. On the
other hand, all the direct summands of $X^{p}$ are perfect and they can have
any object of $\mathsf{P}(\varphi)^{s}$ as JH-factor.

\begin{corollary}
  Any coherent sheaf $\mathcal{F}$ with $\End(\mathcal{F}) = \boldsymbol{k}$
  is stable.
\end{corollary}

\begin{proof}
  The assumption implies that $\mathcal{F}$ is indecomposable.
  If $\mathcal{F}$ were not even semi-stable, it would have at least two
  HN-factors. Using Corollary \ref{cor:sheaves}, we may assume that
  $\varphi_{+}(\mathcal{F})=1$. Thus, $\mathcal{F}$ is a coherent sheaf which
  is neither torsion nor torsion free. This implies that there is a
  non-invertible endomorphism
  $\mathcal{F} \rightarrow \boldsymbol{k}(s) \rightarrow \tors(\mathcal{F})
  \rightarrow \mathcal{F}$, in contradiction to the assumption. Hence,
  $\mathcal{F}\in\mathsf{P}(\varphi)$ is semi-stable. Let
  $\mathcal{J}\in\mathsf{P}(\varphi)$ be its JH-factor.
  From Corollary \ref{cor:sheaves} (\ref{cor:vi}) we obtain a non-zero
  endomorphism $\mathcal{F}\rightarrow \mathcal{J}\rightarrow \mathcal{F}$,
  which can only be an isomorphism, if $\mathcal{F}\cong \mathcal{J}$, so
  $\mathcal{F}$ is indeed stable.  
\end{proof}

The following method can be used to visualise the structure of the category
$\Dbcoh(\boldsymbol{E})$: the vertical slices in Figure \ref{fig:slices} are
thought to correspond to the categories $\mathsf{P}(t)^{s}$ of stable objects. 
  \begin{figure}[hbt]
    \begin{center}
      \setlength{\unitlength}{10mm}
\begin{picture}(11,5)
  \multiput(0,4)(0.2,0){56}{\line(1,0){0.1}}
  \put(0,1){\line(1,0){11.1}}
  \thicklines
  \put(1,1){\line(0,1){3}}\put(1,0.8){\makebox(0,0)[t]{$2$}}
  \put(4,1){\line(0,1){3}}\put(4,0.8){\makebox(0,0)[t]{$1$}}
  \put(7,1){\line(0,1){3}}\put(7,0.8){\makebox(0,0)[t]{$0$}}
  \put(10,1){\line(0,1){3}}\put(10,0.8){\makebox(0,0)[t]{$-1$}}
  \thinlines
  \put(4.9,1){\line(0,1){3}}\put(4.9,0.8){\makebox(0,0)[t]{$t$}}
  \put(1.8,2.5){$\Coh_{\boldsymbol{E}}[1]$}
  \put(5.3,2.5){$\Coh_{\boldsymbol{E}}$}
  \put(7.6,2.5){$\Coh_{\boldsymbol{E}}[-1]$}
\end{picture}
    \end{center}
  \caption{slices}\label{fig:slices}
\end{figure}
They are non-empty if and only if $t\in\mathsf{Q}$, i.e.\/
$\mathbb{R}\exp(\pi it) \cap \mathbb{Z}^{2} \ne \{(0,0)\}$. A point on such a
slice represents a stable object. The extreme stable objects are those which
lie on the dashed upper horizontal line. The labelling below the picture
reflects the phases of the slices. We have chosen to let it decrease from the
left to right in order to have objects with cohomology in negative degrees on
the left and with positive degrees on the right. 

By Proposition \ref{prop:transit}, the group $\widetilde{\SL}(2,\mathbb{Z})$
acts on the set of all stable objects, hence it acts on such pictures.
This action sends slices to slices and acts transitively on the set of 
slices with phase $t\in\mathsf{Q}$. The dashed line of extreme stable objects
is invariant under this action. 

Any indecomposable object $0\ne X\in\Dbcoh(\boldsymbol{E})$ has a
\emph{shadow} in such a picture: it is the set of all stable objects which
occur as JH-factors in the HN-factors of $X$. If this set consists of more
than one point, the shadow is obtained by connecting these points by line
segments.

The following proposition shows that the shadow of an indecomposable object
which consists of more than one point is completely contained in the
extreme line. 
  \begin{figure}[hbt]
    \begin{center}
      \setlength{\unitlength}{10mm}
\begin{picture}(11,5)
  \multiput(0,4)(0.2,0){56}{\line(1,0){0.1}}
  \put(0,1){\line(1,0){11.1}}
  \thicklines
  \put(1,1){\line(0,1){3}}\put(1,0.8){\makebox(0,0)[t]{$2$}}
  \put(4,1){\line(0,1){3}}\put(4,0.8){\makebox(0,0)[t]{$1$}}
  \put(7,1){\line(0,1){3}}\put(7,0.8){\makebox(0,0)[t]{$0$}}
  \put(10,1){\line(0,1){3}}\put(10,0.8){\makebox(0,0)[t]{$-1$}}
  \thinlines
  \put(1.8,2.5){$\Coh_{\boldsymbol{E}}[1]$}
  \put(5.3,2.5){$\Coh_{\boldsymbol{E}}$}
  \put(7.6,2.5){$\Coh_{\boldsymbol{E}}[-1]$}
  \put(4,2){\circle*{0.2}}\put(4.2,2){\makebox(0,0)[l]{$X_{1}$}}
  \put(8.2,3.4){\circle*{0.2}}\put(8.4,3.4){\makebox(0,0)[l]{$X_{2}$}}
  \put(0.3,4){\circle*{0.2}}
  \thicklines\put(0.3,4){\line(1,0){1.5}}
  \put(1.8,4){\circle*{0.2}}
  \thicklines\put(1.8,4){\line(1,0){0.8}}
  \put(2.6,4){\circle*{0.2}}\put(1.3,4.2){\makebox(0,0)[b]{$X_{3}$}}
  \put(4.3,4){\circle*{0.2}}
  \thicklines\put(4.3,4){\line(1,0){1}}
  \put(5.3,4){\circle*{0.2}}\put(4.8,4.2){\makebox(0,0)[b]{$X_{4}$}}
  \put(6.3,4){\circle*{0.2}}\put(6.3,4.2){\makebox(0,0)[b]{$X_{5}$}}
\end{picture}
    \end{center}
  \caption{shadows}\label{fig:example}
\end{figure}

Figure \ref{fig:example} shows the shadows of five different
indecomposable objects:
\begin{itemize}
\item $X_{1}\in\Coh_{\boldsymbol{E}}$ an indecomposable torsion sheaf,
\item $X_{2}\in \Coh_{\boldsymbol{E}}[-1]$ the shift of an indecomposable
  semi-stable locally free sheaf,
\item $X_{3}$ a genuine complex with three extreme HN-factors, one in
  $\Coh_{\boldsymbol{E}}[2]$ and the other two in $\Coh_{\boldsymbol{E}}[1]$,
\item $X_{4}$ an indecomposable torsion free sheaf which is not semi-stable,
\item $X_{5}\in\Coh_{\boldsymbol{E}}$  an indecomposable and semi-stable
  torsion free sheaf which could be perfect or not (a band or a string in the
  language of representation theory).
\end{itemize}
The shadow of an indecomposable object is a single point if and only if this
object is semi-stable.

\begin{proposition}\label{prop:extreme}
  Let $X \in \Dbcoh(\boldsymbol{E})$ be an indecomposable object which is not
  semi-stable. Then, all HN-factors of $X$ are extreme.  
\end{proposition}

\begin{proof}
  Let
  \[\xymatrix@C=.5em{
  0\; \ar[rr] && F_{n}X \ar[rr] \ar[dl]_{\cong}&&   F_{n-1}X \ar[rr]
  \ar[dl]&& \dots \ar[rr] &&F_{1}X   \ar[rr]  &&  F_{0}X \ar@{=}[r]\ar[dl] &
  X\\
  & A_{n} \ar[lu]^{+} && A_{n-1} \ar[lu]^{+} &  & & & & &  A_0 \ar[lu]^{+}&
  }\]
  be a HNF of $X$. If the HN-factor $A_{i}$ were not
  extreme, it could be split into a direct sum $A_{i} \cong A_{i}'
  \oplus A_{i}''$ with $0\ne A_{i}'$ perfect and $A_{i}',
  A_{i}''\in\mathsf{P}(\varphi_{i})$. 
  Because $\varphi_{-}(F_{i+1}X) > \varphi_{i}=\varphi(A_{i}')$, Lemma
  \ref{wesPT:ii} and Serre duality imply
  $$\Hom(A_{i}', F_{i+1}X[1]) \cong \Hom(F_{i+1}X, A_{i}')^{\ast} = 0.$$
  Hence, we can apply Lemma \ref{lem:PengXiao} to the distinguished triangle
  $$F_{i+1}X \rightarrow F_{i}X \rightarrow A_{i}
  \stackrel{+}{\longrightarrow}$$
  and obtain a decomposition $F_{i}X \cong F_{i}'X \oplus A_{i}'$.
  We proceed by descending induction on $j\le i$ to show that there exist
  decompositions $F_{j}X \cong F_{j}'X \oplus A_{i}'$. This is
  obtained from Lemma \ref{lem:PengXiao} applied to the distinguished triangle 
  $$F_{j}'X\oplus A_{i}' \rightarrow F_{j-1}X \rightarrow A_{j-1}
  \stackrel{+}{\longrightarrow}$$
  and using Lemma \ref{wesPT:ii}, Serre duality and
  $\varphi(A_{i}') > \varphi(A_{j-1})$ to get
  $$\Hom(A_{j-1}, A_{i}'[1]) \cong \Hom(A_{i}', A_{j-1})^{\ast} = 0.$$
  We obtain a decomposition $X=F_{0}X \cong F_{0}'X \oplus A_{i}'$ in which
  we have $A_{i}'\ne 0$. Because $X$ was assumed to be indecomposable, we
  should have $X\cong A_{i}'$, but this was excluded by assumption. This
  contradiction shows that all HN-factors $A_{i}$ are necessarily
  extreme. 
\end{proof}

\begin{corollary}\label{cor:types}
  There exist four types of indecomposable objects in the category
  $\Coh_{\boldsymbol{E}}$: 
  \begin{enumerate}
  \item \label{type:i} semi-stable with perfect JH-factor;
  \item \label{type:ii} semi-stable, perfect but its JH-factor extreme;
  \item \label{type:iii} semi-stable and extreme;
  \item \label{type:iv} not semi-stable, with all its HN-factors extreme.
  \end{enumerate}
\end{corollary}

A similar statement is true for $\Dbcoh(\boldsymbol{E})$. In this case, the
objects of types (\ref{type:i}), (\ref{type:ii}) and (\ref{type:iii}) are
shifts of coherent sheaves, whereas genuine complexes are possible for objects
of type (\ref{type:iv}). 
Types (\ref{type:ii}), (\ref{type:iii}) and (\ref{type:iv}) were not available
in the smooth case. 

Examples of type (\ref{type:i}) are simple vector bundles and structure sheaves
$\boldsymbol{k}(x)$ of smooth points $x\in\boldsymbol{E}$. All indecomposable
objects with a shadow not on the extreme line fall into type (\ref{type:i}).
Under the equivalences of Corollary \ref{cor:equiv}, indecomposable semi-stable
locally free sheaves with extreme JH-factor correspond, in the nodal case,
precisely to those torsion sheaves with support at the node $s$, which are
called bands, see \cite{BurbanKreussler}. 
Examples of type (\ref{type:iii}) are the stable coherent sheaves which are not
locally free and the structure sheaf $\boldsymbol{k}(s)$ of the singular point
$s\in\boldsymbol{E}$. Moreover, in the nodal case, the torsion sheaves with
support at $s$, which are called strings in \cite{BurbanKreussler}, are of
type (\ref{type:iii}) as well. Examples of objects of type (\ref{type:iv}) are
given below.  

\begin{example}
  We shall construct torsion free sheaves on nodal $\boldsymbol{E}$ with an
  arbitrary finite number of HN-factors. This implies that the number of points
  in a shadow of an indecomposable object in $\Dbcoh(\boldsymbol{E})$ is not
  bounded.

  Recall from \cite{DrozdGreuel} that any indecomposable torsion free sheaf
  which is not locally free, is isomorphic to a sheaf
  $\mathcal{S}(\boldsymbol{d}) = p_{n\ast} \mathcal{L}(\boldsymbol{d})$. We use
  here the notation of \cite{BurbanKreussler}, Section 3.5, so that $p_{n}:
  \boldsymbol{I_{n}} \rightarrow \boldsymbol{E}$ denotes a certain morphism
  from the chain $\boldsymbol{I_{n}}$ of $n$ smooth rational curves to the
  nodal curve $\boldsymbol{E}$. If $\boldsymbol{d}=(d_{1},\ldots,d_{n})
  \in\mathbb{Z}^{n}$, we denote by $\mathcal{L}(\boldsymbol{d})$ the line
  bundle on $\boldsymbol{I_{n}}$ which has degree $d_{\nu}$ on the $\nu$-th
  component of $\boldsymbol{I_{n}}$. We know $\rk(\mathcal{S}(\boldsymbol{d}))
  = n$ and $\deg(\mathcal{S}(\boldsymbol{d})) = 1+\sum d_{\nu}$. We obtain,
  in particular, that for any $\varphi\in\mathsf{Q} \cap (0,1)$ there exist
  $n\in\mathbb{Z}$ and $\boldsymbol{d}(\varphi)\in\mathbb{Z}^{n}$ such that
  $\mathcal{S}(\boldsymbol{d}(\varphi))$ is the unique extreme element in
  $\mathsf{P}(\varphi)^{s}$. On the other hand, if $\boldsymbol{d}'\in
  \mathbb{Z}^{n'}, \boldsymbol{d}''\in \mathbb{Z}^{n''}$ and
  $\boldsymbol{d} = (\boldsymbol{d}_{+}', \boldsymbol{d}'')\in
  \mathbb{Z}^{n'+n''}$, where $\boldsymbol{d}_{+}'$ is obtained from
  $\boldsymbol{d}'$ by adding $1$ to the last component, we have an exact
  sequence
  $$0\rightarrow \mathcal{S}(\boldsymbol{d}') \rightarrow
  \mathcal{S}(\boldsymbol{d}) \rightarrow
  \mathcal{S}(\boldsymbol{d}'') \rightarrow 0$$
  see for example \cite{Mozgovoy}. Hence, if we start with a sequence
  $0<\varphi_{0} <\varphi_{1}< \ldots <\varphi_{m} <1$ where
  $\varphi_{\nu}\in\mathsf{Q}$ and define 
  $$\boldsymbol{d}^{(m)} = \boldsymbol{d}(\varphi_{m})\quad\text{ and }\quad
  \boldsymbol{d}^{(\nu)} =
  (\boldsymbol{d}_{+}^{(\nu+1)},\boldsymbol{d}(\varphi_{\nu})) \text{ for }
  m > \nu \ge 0,$$
  we obtain an indecomposable torsion free sheaf
  $\mathcal{S}(\boldsymbol{d}^{(0)})$ whose HN-factors are the extreme stable
  sheaves $\mathcal{S}(\boldsymbol{d}(\varphi_{\nu})) \in
  \mathsf{P}(\varphi_{\nu}), 0\le \nu \le m$. The HNF of this sheaf is given by
  $$\mathcal{S}(\boldsymbol{d}^{(m)}) \subset
  \mathcal{S}(\boldsymbol{d}^{(m-1)}) \subset \ldots \subset
  \mathcal{S}(\boldsymbol{d}^{(0)}).$$
  The sheaf $\mathcal{S}(\boldsymbol{d}^{(0)})$ is of type (\ref{type:iv}) and
  not perfect. 
\end{example}

\begin{example}\label{ex:extremefactors}
  Suppose $\boldsymbol{E}$ is nodal and let
  $\pi:C_{2}\rightarrow\boldsymbol{E}$ be an \'etale morphism of degree 
  two, where $C_{2}$ denotes a reducible curve which has two components, both
  isomorphic to $\mathbb{P}^{1}$ and which intersect transversally at two 
  distinct points.
  By $i_{\nu}:\mathbb{P}^{1}\rightarrow \boldsymbol{E},\;\nu=1,2$ we denote the
  morphisms which are induced by the embeddings of the two components of
  $C_{2}$. 
  There is a $\boldsymbol{k}^{\times}$-family of
  line bundles on $C_{2}$, whose restriction to one component is
  $\mathcal{O}_{\mathbb{P}^{1}}(-2)$ and to the other is
  $\mathcal{O}_{\mathbb{P}^{1}}(2)$. The element in $\boldsymbol{k}^{\times}$ 
  corresponds to a gluing parameter over one of the two singularities of
  $C_{2}$. If $\mathcal{L}$ denotes one such line bundle,
  $\mathcal{E}:=\pi_{\ast}\mathcal{L}$ is an 
  indecomposable vector bundle of rank two and degree zero on $\boldsymbol{E}$.
  Let us fix notation so that $i_{1}^{\ast}\mathcal{E}
  \cong \mathcal{O}_{\mathbb{P}^{1}}(-2)$ and $i_{2}^{\ast}\mathcal{E} \cong
  \mathcal{O}_{\mathbb{P}^{1}}(2)$. There is an exact sequence of coherent
  sheaves on $\boldsymbol{E}$ 
  \begin{equation}\label{eq:nonssvb}
  0\rightarrow i_{2\ast} \mathcal{O}_{\mathbb{P}^{1}} \rightarrow
  \mathcal{E} 
  \rightarrow i_{1\ast} \mathcal{O}_{\mathbb{P}^{1}}(-2) \rightarrow 0.
  \end{equation}
  Because the torsion free sheaves $i_{2\ast} \mathcal{O}_{\mathbb{P}^{1}}$ and
  $i_{1\ast} \mathcal{O}_{\mathbb{P}^{1}}(-2)$ have rank one and
  $\boldsymbol{E}$ is irreducible, they are stable. Because $\varphi(i_{2\ast}
  \mathcal{O}_{\mathbb{P}^{1}}) = 3/4$ and $\varphi(i_{1\ast}
  \mathcal{O}_{\mathbb{P}^{1}}(-2)) = 1/4$, Theorem \ref{thm:uniqueHNF}
  implies that the HNF of $\mathcal{E}$ is given by the 
  exact sequence (\ref{eq:nonssvb}). The HN-factors are the two torsion
  free sheaves of rank one $i_{2\ast} \mathcal{O}_{\mathbb{P}^{1}}$ and
  $i_{1\ast} \mathcal{O}_{\mathbb{P}^{1}}(-2)$, which are not locally
  free. These are the extreme stable elements with phases $3/4$ and $1/4$
  respectively. Therefore, the indecomposable vector bundle $\mathcal{E}$ is a
  perfect object of type (\ref{type:iv}) which satisfies
  $\varphi_{-}(\mathcal{E})=1/4$  and $\varphi_{+}(\mathcal{E})=3/4$.
\end{example}

\begin{remark}\label{rem:notperfect}
  This example shows that the full sub-category of perfect complexes in the
  category $\Dbcoh(\boldsymbol{E})$ is not closed under taking
  Harder-Narasimhan factors. We interpret this to be an indication that the
  derived category of perfect complexes is not an appropriate object for
  homological mirror symmetry on singular Calabi-Yau varieties.
\end{remark}

\begin{remark}
  It seems plausible that methods similar to those of this section could be
  applied to study the derived category of representations of certain derived
  tame associative algebras. Such may include gentle algebras, skew-gentle
  algebras and degenerated tubular algebras. 
  The study of Harder-Narasimhan filtrations in conjunction with the action of
  the group of exact auto-equivalences of the derived category may provide new
  insight into the combinatorics of indecomposable objects in these derived
  categories.
\end{remark}

\begin{proposition}\label{wesPT}
  Suppose $X,Y\in \Dbcoh(\boldsymbol{E})$ are non-zero.
  \begin{enumerate}
  \item \label{wesPT:iii} If $\varphi_{-}(X) <
    \varphi_{+}(Y) < \varphi_{-}(X)+1$, then $\Hom(X,Y)\ne 0$.
  \item \label{wesPT:iv} If $X$ and $Y$ are indecomposable objects which are
    not of type (\ref{type:i}) in Corollary \ref{cor:types} and which satisfy
    $\varphi_{-}(X) = \varphi_{+}(Y)$, then $\Hom(X,Y)\ne 0$.  
  \end{enumerate}
\end{proposition}

\begin{proof}
  If $X$ and $Y$ are semi-stable objects, the claim
  (\ref{wesPT:iii}) was proved in Corollary \ref{cor:sheaves} (\ref{cor:iii}).
  Similarly, (\ref{wesPT:iv}) for two semi-stable objects follows from
  Corollary \ref{cor:sheaves} (\ref{cor:vi}), because there is only one
  non-perfect object in $\mathsf{P}(\varphi)^{s}$.

  For the rest of the proof we treat both cases,  (\ref{wesPT:iii}) and
  (\ref{wesPT:iv}) simultaneously.
  For the proof of (\ref{wesPT:iv}) we keep in mind that Proposition
  \ref{prop:extreme} implies that no HN-factor has a perfect summand, if
  the object is indecomposable but not semi-stable.
  If $X\in\mathsf{P(\varphi)}$ is semi-stable but $Y\in\Dbcoh(\boldsymbol{E})$
  is arbitrary, we let 
    \[\xymatrix@C=.5em{
    0\; \ar[rr] && F_{m}Y \ar[rr] \ar[dl]_{\cong}&&   F_{m-1}Y \ar[rr]
    \ar[dl]&& \dots \ar[rr] &&F_{1}Y   \ar[rr]  &&  F_{0}Y \ar@{=}[r]\ar[dl] &
    Y\\
    & B_{m} \ar[lu]^{+} && B_{m-1} \ar[lu]^{+} &  & & & & &  B_0 \ar[lu]^{+}&
  }\]
  be a HNF of $Y$. As $\varphi(B_{m})=\varphi_{+}(Y)$ we know already
  $\Hom(X,B_{m})\ne 0$. 
  By assumption, we have $\varphi(B_{i}[-1]) = \varphi(B_{i})
  -1 \le \varphi_{+}(Y)-1 < \varphi(X)$. Hence, by Lemma \ref{wesPT:ii},
  $\Hom(X, B_{i}[-1]) =0$ and the cohomology sequence of the distinguished
  triangle $F_{i+1}Y\rightarrow F_{i}Y\rightarrow B_{i}
  \stackrel{+}{\rightarrow}$ provides an inclusion $\Hom(X, F_{i+1}Y) \subset
  \Hom(X, F_{i}Y)$. This implies $0\ne \Hom(X,B_{m})\subset \Hom(X,Y)$.

  Finally, in the general case, we let
    \[\xymatrix@C=.5em{
    0\; \ar[rr] && F_{n}X \ar[rr] \ar[dl]_{\cong}&&   F_{n-1}X \ar[rr]
    \ar[dl]&& \dots \ar[rr] &&F_{1}X   \ar[rr]  &&  F_{0}X \ar@{=}[r]\ar[dl] &
    X\\
    & A_{n} \ar[lu]^{+} && A_{n-1} \ar[lu]^{+} &  & & & & &  A_0 \ar[lu]^{+}&
  }\]
  be a HNF of $X$. As $\varphi(A_{0})=\varphi_{-}(X)$ we have
  $\Hom(A_{0},Y)\ne 0$.
  Because $\varphi_{-}(F_{1}X[1]) = \varphi_{-}(F_{1}X) +1 =
  \varphi(A_{1}) +1 > \varphi_{-}(X)+1 > \varphi_{+}(Y)$, Lemma
  \ref{wesPT:ii} implies $\Hom(F_{1}X[1], Y)=0$. The distinguished triangle
  $F_{1}X \rightarrow X \rightarrow A_{0} \stackrel{+}{\rightarrow}$ gives us
  now an inclusion $0\ne \Hom(A_{0},Y) \subset \Hom(X,Y)$ and so the claim.
\end{proof}

In \cite{YangBaxter}, Polishchuk asked for the classification of all
spherical objects in the bounded derived category of a singular projective
curve of arithmetic genus one. Below, we shall solve this problem for
irreducible curves.

Let $\boldsymbol{E}$ be an irreducible projective curve of arithmetic genus
one over our base field $\boldsymbol{k}$. 
Recall that in this  case an object $X\in\Dbcoh(\boldsymbol{E})$ is
\emph{spherical} if 
$$X \text{ is perfect and }\quad
\Hom(X,X[i]) \cong
\begin{cases}
  \boldsymbol{k} & \text{if }\;  i \in     \{0,1\} \\
          0      & \text{if }\;  i \not\in  \{0,1\} 
\end{cases}
$$

\begin{proposition}\label{prop:spherical}
  Let $\boldsymbol{E}$ be an irreducible projective curve of arithmetic genus
  one  and $X\in\Dbcoh(\boldsymbol{E})$. Then the following  are equivalent:
  \begin{enumerate}
  \item\label{spher:i} $X$ is spherical;
  \item \label{spher:ii}$\Hom(X,X[i]) \cong
    \begin{cases}
      \boldsymbol{k} & \text{if }\;  i = 0 \\
              0      & \text{if }\;  i = 2 \;\text{ or }\; i<0;
    \end{cases}$
  \item\label{spher:iii} $X$ is perfect and stable;
  \item\label{spher:iv} there exists $n\in\mathbb{Z}$ such that $X[n]$ is
  isomorphic to a simple vector bundle or to a torsion sheaf of length one
  which is supported at a smooth point of $\boldsymbol{E}$. 
  \end{enumerate}
  In particular, the group of exact auto-equivalences of
  $\Dbcoh(\boldsymbol{E})$ acts transitively on the set of all spherical
  objects. 
\end{proposition}

\begin{proof}
  The implication (\ref{spher:i})$\Rightarrow$(\ref{spher:ii}) is obvious.

  Let us prove (\ref{spher:ii})$\Rightarrow$(\ref{spher:iii}).
  First, we observe that $\Hom(X,X) \cong \boldsymbol{k}$ implies that
  $X$ is indecomposable. Suppose, $X$ is not semi-stable. This is equivalent
  to $\varphi_{+}(X)>\varphi_{-}(X)$. By Proposition \ref{prop:extreme}
  we know that all HN-factors of $X$ are extreme.
  Let $M\ge 0$ be the unique integer with $M\le \varphi_{+}(X) -
  \varphi_{-}(X) < M+1$.

  If $M< \varphi_{+}(X) - \varphi_{-}(X) <M+1$,
  Proposition \ref{wesPT} (\ref{wesPT:iii}) implies $\Hom(X,X[-M])\ne
  0$. Under the assumption (\ref{spher:ii}), this is possible only if $M=0$.
  On the other hand, if $M=\varphi_{+}(X) - \varphi_{-}(X)$, we obtain from
  Proposition \ref{wesPT} (\ref{wesPT:iv}) $\Hom(X,X[-M]) \ne 0$.
  Again, this implies $M=0$.
  So, we have $0< \varphi_{+}(X) - \varphi_{-}(X) <1$.

  If we apply the functor $\Hom(\ARG,X)$ to
  $F_{1}X\stackrel{u}{\rightarrow} X \rightarrow A_{0}
  \stackrel{+}{\longrightarrow}$,
  the rightmost distinguished triangle of the HNF of $X$, we obtain the exact
  sequence 
  $$\Hom(F_{1}X[1],X) \rightarrow \Hom(A_{0},X) \rightarrow \Hom(X,X)
  \rightarrow \Hom(F_{1}X,X),$$
  in which the leftmost term $\Hom(F_{1}X[1],X)=0$ by Lemma \ref{wesPT:ii},
  because $\varphi_{-}(F_{1}X[1]) > \varphi_{-}(X)+1>\varphi_{+}(X)$. The third
  morphism in this sequence is not the zero map, as it sends $\mathsf{Id}_{X}$
  to $u\ne 0$. Because $\Hom(X,X)$ is one dimensional, this is only possible
  if $\Hom(A_{0},X)=0$. But Proposition \ref{wesPT} (\ref{wesPT:iii}) and
  $\varphi(A_{0})<\varphi_{+}(X)< \varphi(A_{0})+1$ imply $\Hom(A_{0},X)\ne 0$.
  This contradiction shows that $X$ must be semi-stable.
  
  We observed earlier that all the
  JH-factors of an indecomposable semi-stable object are isomorphic to each
  other. Therefore, any indecomposable semi-stable object which is not stable
  has a space of endomorphisms of dimension at least two. So, we conclude 
  $X\in\mathsf{P}(\varphi)^{s}$ for some $\varphi\in\mathbb{R}$.
  
  Because $\Hom(\boldsymbol{k}(s),\boldsymbol{k}(s)[2]) \cong
  \Ext^{2}(\boldsymbol{k}(s),\boldsymbol{k}(s))\ne 0$, the transitivity of the
  action of $\widetilde{\SL}(2,\mathbb{Z})$ on the set $\mathsf{Q}$ implies
  that none of the extreme stable objects satisfies the condition
  (\ref{spher:ii}). Hence, $X$ is perfect and stable. 

  To prove (\ref{spher:iii})$\Rightarrow$(\ref{spher:i}), we observe that
  the group of automorphisms of the curve $\boldsymbol{E}$ acts
  transitively on the regular locus $\boldsymbol{E}\setminus\{s\}$. Hence, by
  Proposition \ref{prop:transit}, the group of auto-equivalences of
  $\Dbcoh(\boldsymbol{E})$ acts transitively on the set of all perfect stable
  objects. Because, for example, the structure sheaf
  $\mathcal{O}_{\boldsymbol{E}}$ is spherical, it is now clear that all
  perfect stable objects are indeed spherical and that the group of 
  exact auto-equivalences of $\Dbcoh(\boldsymbol{E})$ acts transitively on the
  set of all spherical objects.

  To show the equivalence with (\ref{spher:iv}), it remains to recall that any
  perfect coherent torsion free sheaf on $\boldsymbol{E}$ is locally free. This
  follows easily from the Auslander-Buchsbaum formula because we are working in
  dimension one. 
\end{proof}

\section{Description of $t$-structures in the case of a 
  singular Weierstra\ss{} curve}\label{sec:tstruc}

The main result of this section is a description of all $t$-structures on the
derived category of a singular Weierstra\ss{} curve $\boldsymbol{E}$. This
generalises results of \cite{GRK} and  \cite{Pol1}, where the smooth case was
studied. As an application, we obtain a description of the group
$\Aut(\Dbcoh(\boldsymbol{E}))$ of all exact auto-equivalences of
$\Dbcoh(\boldsymbol{E})$. A second application is a description of Bridgeland's
space of stability conditions on $\boldsymbol{E}$.

Recall that a $t$-structure on  a triangulated category $\mathsf{D}$ is a pair
of full subcategories $(\mathsf{D}^{\le 0}, \mathsf{D}^{\ge 0})$ such that,
with the notation $\mathsf{D}^{\ge n} := \mathsf{D}^{\ge
0}[-n]$ and $\mathsf{D}^{\le n} := \mathsf{D}^{\le 0}[-n]$ for any $n\in
\mathbb{Z}$, the following holds:
\begin{enumerate}
\item $\mathsf{D}^{\le 0} \subset \mathsf{D}^{\le 1}$ and $\mathsf{D}^{\ge 1}
  \subset   \mathsf{D}^{\ge 0}$;
\item $\Hom(\mathsf{D}^{\le 0}, \mathsf{D}^{\ge 1}) = 0$;
\item\label{def:tiii} for any object $X \in \mathsf{D}$ there exists a
  distinguished triangle
  $$A \rightarrow X \rightarrow B \stackrel{+}{\longrightarrow}$$
  with $A \in \mathsf{D}^{\le 0}$ and $B \in \mathsf{D}^{\ge 1}.$
\end{enumerate}

If $(\mathsf{D}^{\le 0}, \mathsf{D}^{\ge 0})$ is a $t$-structure then ${\sf
A} = \mathsf{D}^{\le 0} \cap \mathsf{D}^{\ge 0}$ has a structure of an
Abelian category. It is called the \emph{heart} of the $t$-structure. 
In this way, $t$-structures on the derived category
$\Dbcoh(\boldsymbol{E})$ lead to interesting Abelian categories embedded into
it. The natural $t$-structure on $\Dbcoh(\boldsymbol{E})$ has $\mathsf{D}^{\le
n}$ equal to the full subcategory formed by all complexes with non-zero
cohomology in degree less or equal to $n$ only. Similarly, the full subcategory
$\mathsf{D}^{\ge n}$ consists of all complexes $X$ with $H^{i}(X)=0$ for all
$i<n$. The heart of the natural $t$-structure is the Abelian category
$\Coh_{\boldsymbol{E}}$.

In addition to the natural $t$-structure we also have many interesting
$t$-structures on $\Dbcoh(\boldsymbol{E})$. 
In order to describe them, we introduce the following notation. We continue to
work with the notion of stability and the notation introduced in the previous
section. 
If $\mathsf{P}\subset\mathsf{P}(\theta)^{s}$ is a subset, we denote by
$\mathsf{D}[\mathsf{P}, \infty)$ the full subcategory of
$\Dbcoh(\boldsymbol{E})$ which is defined as follows:
$X\in\Dbcoh(\boldsymbol{E})$ is in $\mathsf{D}[\mathsf{P},
\infty)$ if and only if $X=0$ or all its HN-factors, which have at least one
JH-factor which is not in $\mathsf{P}$, have phase $\varphi>\theta$.
Similarly, $\mathsf{D}(-\infty,\mathsf{P}]$ denotes the category which is
generated by $\mathsf{P}$ and all $\mathsf{P}(\varphi)$ with
$\varphi<\theta$. 
If $\mathsf{P}=\mathsf{P}(\theta)^{s}$ we may abbreviate
$\mathsf{D}[\theta,\infty) = \mathsf{D}[\mathsf{P}, \infty)$ and
$\mathsf{D}(-\infty,\theta] = \mathsf{D}(-\infty,\mathsf{P}]$. 
Similarly, if $\mathsf{P}=\emptyset$ we use
the abbreviations $\mathsf{D}(\theta,\infty)$ and $\mathsf{D}(-\infty,\theta)$.
For any open, closed or half-closed interval $I\subset\mathbb{R}$ we define
the full subcategories $\mathsf{D}I$ precisely in the same way. Thus, an
object $0\ne X\in\Dbcoh(\boldsymbol{E})$ is in $\mathsf{D}I$ if and only if
$\varphi_{-}(X)\in I$ and $\varphi_{+}(X)\in I$.

\begin{proposition}\label{prop:texpl}
  Let $\theta\in\mathbb{R}$ and $\mathsf{P}(\theta)^{-} \subset
  \mathsf{P}(\theta)^{s}$ be arbitrary. Denote by $\mathsf{P}(\theta)^{+} =
  \mathsf{P}(\theta)^{s} \setminus \mathsf{P}(\theta)^{-}$ the complement of
  $\mathsf{P}(\theta)^{-}$. Then,  
  $$\mathsf{D}^{\le0} :=  \mathsf{D}[\mathsf{P}(\theta)^{-}, \infty)$$
  defines a $t$-structure on $\Dbcoh(\boldsymbol{E})$ with
  $$\mathsf{D}^{\ge1} :=  \mathsf{D}(-\infty,\mathsf{P}(\theta)^{+}].$$
  The heart $\mathsf{A}(\theta,\mathsf{P}(\theta)^{-})$ of it is
  the category $\mathsf{D}[\mathsf{P}(\theta)^{-},
  \mathsf{P}(\theta)^{+}[1]]$, which consists of those objects
  $X\in\Dbcoh(\boldsymbol{E})$ whose HN-factors either have
  phase $\varphi\in(\theta,\theta+1)$ or have all its JH-factors in
  $\mathsf{P}(\theta)^{-}$ or $\mathsf{P}(\theta)^{+}[1]$.
\end{proposition}

\begin{proof}
  The only non-trivial property which deserves a proof is (\ref{def:tiii}) in
  the definition of $t$-structure. Given $X\in \Dbcoh(\boldsymbol{E})$, we
  have to show that there exists a distinguished triangle $A\rightarrow X
  \rightarrow B \stackrel{+}{\rightarrow}$ with $A\in \mathsf{D}^{\le0}$ and
  $B\in \mathsf{D}^{\ge1}$. In order to construct it, let
  \[\xymatrix@C=.5em{
    0\; \ar[rr] && F_{n}X \ar[rr] \ar[dl]_{\cong}&&   F_{n-1}X \ar[rr]
    \ar[dl]&& \dots \ar[rr] &&F_{1}X   \ar[rr]  &&  F_{0}X \ar@{=}[r]\ar[dl] &
    X\\
    & A_{n} \ar[lu]^{+} && A_{n-1} \ar[lu]^{+} &  & & & & &  A_0 \ar[lu]^{+}&
  }\]
  be the HNF of $X$. Because $\varphi(A_{i+1})>\varphi(A_{i})$ for all $i$,
  there exists an integer $k$, $0\le k \le n+1$ such that $\varphi(A_{k})\ge
  \theta >\varphi(A_{k-1})$. If $\varphi(A_{k})>\theta$, this implies
  $A_{i}\in \mathsf{D}^{\le0}$, if $i\ge k$ and $A_{i}\in \mathsf{D}^{\ge1}$,
  if $i<k$. In particular, $F_{k}X\in \mathsf{D}^{\le0}$. In this case, we
  define $A:=F_{k}X$ and let $A=F_{k}X \rightarrow X$ be the composition of
  the morphisms in the HNF. 
  If, however, $\varphi(A_{k})=\theta$, there is a splitting $A_{k}\cong
  A_{k}^{-} \oplus A_{k}^{+}$ such that all JH-factors of $A_{k}^{-}$ (resp.\/
  $A_{k}^{+}$) are in $\mathsf{P}(\theta)^{-}$ (resp.\/
  $\mathsf{P}(\theta)^{+}$).
  Now, we apply Lemma \ref{lem:connect} to the distinguished triangles
  $F_{k+1}X \stackrel{f}{\longrightarrow} F_{k}X \longrightarrow A_{k}
  \stackrel{+}{\longrightarrow}$
  and
  $A_{k}^{-} \longrightarrow A_{k} \longrightarrow A_{k}^{+}
  \stackrel{+}{\longrightarrow}$, given by the splitting of $A_{k}$, 
  to obtain a factorisation $F_{k+1}X \rightarrow A \rightarrow F_{k}X$ of
  $f$ and two distinguished triangles  
  $$\xymatrix@C=.5em{F_{k+1}X \ar[rr] && A \ar[dl]\ar[rr] && F_{k}X.\ar[dl]\\
    & A_{k}^{-} \ar[ul]^{+} && A_{k}^{+}\ar[ul]^{+}}$$
  Part of the given HNF of $X$ together with the left one of these two
  triangles form a  HNF of $A$, whence $A\in \mathsf{D}^{\le0}$.
  Again, we let $A\rightarrow X$
  be obtained by composition with the morphisms in the HNF of $X$. 
  In any case, we choose a distinguished triangle $A\rightarrow X \rightarrow B
  \stackrel{+}{\rightarrow}$, where $A\rightarrow X$ is the morphism chosen
  before. From Lemma \ref{lem:split} or Remark \ref{rem:split} we obtain $B\in
  \mathsf{D}^{\ge1}$. This proves the proposition.
\end{proof}
We shall also need the following standard result. 

\begin{lemma}\label{lem:tsummands}
  Let $(\mathsf{D}^{\le 0}, \sf{D}^{\ge 0})$ be a $t$-structure on a
  triangulated category. If $X \oplus Y \in \sf{D}^{\le 0}$ then
  $X \in \sf{D}^{\le 0}$ and $Y \in \sf{D}^{\le 0}$.
  The corresponding statement holds for $\sf{D}^{\ge 0}$.
\end{lemma}

\begin{proof}
  Let $A \stackrel{f}{\longrightarrow} X \stackrel{g}{\longrightarrow} B
  \stackrel{+}{\longrightarrow}$ be a distinguished triangle with
  $A\in\mathsf{D}^{\le 0}$ and $B\in\mathsf{D}^{\ge 1}$, which exists due to
  the definition of a $t$-structure. If $X\not\in \mathsf{D}^{\le 0}$, we
  necessarily have $g\ne 0$ and $B \ne 0$.
  Because $\Hom(\sf{D}^{\le 0}, \sf{D}^{\ge 1}) = 0$, the composition $X
  \oplus Y \stackrel{p}{\longrightarrow} X \stackrel{g}{\longrightarrow} B$,
  in which $p$ denotes the natural projection,  must be zero. If
  $i:X\rightarrow X\oplus Y$ denotes the canonical morphism, we 
  obtain $g=g\circ p\circ i=0$, a contradiction. In the same way it follows
  that $Y\in\mathsf{D}^{\le 0}$.
\end{proof}

Recall that an Abelian category is called \emph{Noetherian}, if any sequence of
epimorphisms stabilises, this means that for any sequence of epimorphisms
$f_{k}:A_{k}\rightarrow A_{k+1}$ there exists an integer $k_{0}$ such that
$f_{k}$ is an isomorphism for all $k\ge k_{0}$. 

\begin{lemma}\label{lem:heart}
  The heart $\mathsf{A}(\theta,\mathsf{P}(\theta)^{-})$ of the $t$-structure,
  which was described in Proposition \ref{prop:texpl}, is Noetherian if and
  only if $\mathsf{P}(\theta) \ne \{0\}$ and $\mathsf{P}(\theta)^{-} =
  \emptyset$. In this case, ${\sf A}(\theta, \emptyset) = \mathsf{D}(\theta,
  \theta +1]$.  
\end{lemma}

\begin{proof}  
  If $\mathsf{P}(\theta) = \{0\}$ then
  $\mathsf{A}(\theta,\mathsf{P}(\theta)^{-}) = 
  \mathsf{D}(\theta,\theta+1)$. This category is not Noetherian.
  To prove this, we follow the proof of Polishchuk in the smooth case
  \cite{Pol1}, Proposition 3.1.
  We are going to show for any non-zero locally free shifted sheaf $E \in
  \mathsf{D}(\theta, \theta +1)$, the 
  existence of a locally free shifted sheaf $F$ and an epimorphism
  $E\twoheadrightarrow F$ in $\mathsf{D}(\theta, \theta +1)$, which is not an
  isomorphism. This will be sufficient to show that $\mathsf{D}(\theta, \theta
  +1)$ is not Noetherian.

  By applying an appropriate shift, we may assume $0<\theta<1$. Under this
  assumption, for every stable coherent sheaf $G$ we have
  \begin{align*}
    G\in\mathsf{D}(\theta, \theta +1) &\iff \theta<\varphi(G)\le 1\\
    G[1]\in\mathsf{D}(\theta, \theta +1) &\iff 0<\varphi(G)< \theta.
  \end{align*}
  For any two objects $X,Y\in\Dbcoh(\boldsymbol{E})$ we define the Euler form
  to be
  $$\langle X,Y \rangle = \rk(X)\deg(Y) - \deg(X)\rk(Y)$$
  which is the imaginary part of $\overline{Z(X)}Z(Y)$. If $X$ and $Y$ are
  coherent sheaves and one of them is perfect, we have
  $$\langle X,Y \rangle = \chi(X,Y) := \dim\Hom(X,Y) - \dim \Ext^{1}(X,Y).$$
  This remains true, if we apply arbitrary shifts to the sheaves $X,Y$, where
  we understand $\chi(X,Y)=\sum_{\nu} (-1)^{\nu} \dim \Hom(X,Y[\nu]).$

  Let $E\in\mathsf{D}(\theta, \theta +1)$ be an arbitrary non-zero locally
  free shifted sheaf. We look at the strip in the plane between the lines
  $L(0):= \mathbb{R}\exp(i\pi\theta)$ and $L(E):=L(0)+Z(E)$. This strip must
  contain lattice points in its interior. 
    \begin{figure}[hbt]
      \begin{center}
        \setlength{\unitlength}{10mm}
  \begin{picture}(11,6)
    \put(1.5,2){\vector(1,0){9.5}}\put(11,1.9){\makebox(0,0)[t]{$-\deg$}}
    \put(6,0){\vector(0,1){6}}\put(5.8,6){\makebox(0,0)[r]{$\rk$}}
    \put(2,0){\line(2,1){9}}
    \put(11,4.3){\makebox(0,0)[t]{$\theta$}}
    \put(2.8,0){\makebox(0,0)[l]{$\theta+1$}}
    \put(1,1.5){\line(2,1){9}}
    \put(6,2){\vector(-1,1){1}}\put(4.8,2.9){\makebox(0,0)[t]{$F$}}
    \put(6,2){\vector(2,3){2}}\put(8.2,4.9){\makebox(0,0)[t]{$E$}}
    \put(5,3){\vector(3,2){3}}
    \put(10.5,4.5){\makebox(0,0)[b]{$L(0)$}}
    \put(10.5,6){\makebox(0,0)[t]{$L(E)$}}
  \end{picture}
      \end{center}
    \caption{}\label{fig:strip}\end{figure}
  Therefore, there exists a lattice point $Z_{F}$ in this strip which enjoys
  the following properties:
  \begin{enumerate}
  \item\label{nopoint} the only lattice points on the closed triangle whose
  vertices are $0, Z(E), Z_{F}$, are its vertices;
  \item\label{phase} $\varphi_{F} > \varphi(E)$.
  \end{enumerate}
  By $\varphi_{F}$ we denote here the unique number which satisfies
  $\theta <\varphi_{F} < \theta+1$ and $Z_{F}\in \mathbb{R}\exp(i\pi
  \varphi_{F})$.
  Because $\SL(2,\mathbb{Z})$ acts transitively on $\mathsf{Q}$, there exists
  a stable non-zero locally free shifted sheaf $F\in\mathsf{D}(\theta, \theta
  +1)$ with $Z(F)=Z_{F}$ and $\varphi(F)=\varphi_{F}$.
  The assumption $\mathsf{P}(\theta)=\{0\}$ implies
  $\mathbb{R}\exp(i\pi\theta)\cap\mathbb{Z}^{2} = \{0\}$, hence, $Z(E)$ is the
  only lattice point on the line $L(E)$. This implies that $Z(F)$ is not on
  the boundary of the stripe between $L(0)$ and $L(E)$. In particular,
  $Z(E)-Z(F)$ is contained in the same half-plane of $L(0)$ as $Z(E)$ and
  $Z(F)$, see Figure \ref{fig:strip}.
  Condition (\ref{nopoint}) implies $\langle E,F \rangle = 1$. Because $E$ is
  locally free, condition (\ref{phase}) implies
  $$\Ext^{1}(E,F) = \Hom(F,E) = 0.$$
  Hence, $\Hom(E,F)\cong \boldsymbol{k}$. The evaluation map gives, therefore,
  a distinguished triangle
  $$\Hom(E,F)\otimes E \rightarrow F \rightarrow T_{E}(F)
  \stackrel{+}{\longrightarrow}$$
  with $T_{E}(F)\in \Dbcoh(\boldsymbol{E})$.
  If $C:=T_{E}(F)[-1]$ we obtain a distinguished triangle
  \begin{equation}
    \label{eq:mutation}
    C\rightarrow E\rightarrow F \stackrel{+}{\longrightarrow}
  \end{equation}
  with $Z(C)=Z(E)-Z(F)$.
  Because $E$ is a stable non-zero shifted locally free sheaf, it is spherical
  by Proposition \ref{prop:spherical} and so $T_{E}$ is an equivalence. This
  implies that $T_{E}(F)$ is spherical and, by Proposition
  \ref{prop:spherical} again, $C$ is a stable non-zero shifted locally free
  sheaf. All morphisms in the distinguished triangle (\ref{eq:mutation}) are
  non-zero because $C, E, F$ are indecomposable, see Lemma \ref{lem:PengXiao}.
  Using Lemma \ref{wesPT:ii}, this implies $\theta-1<\varphi(C)<\theta+1$.
  However, we have seen in which half-plane $Z(C)$ is contained, so that we
  must have $\theta<\varphi(C)<\theta+1$, which implies
  $C\in\mathsf{D}(\theta,\theta+1)$. The distinguished triangle
  (\ref{eq:mutation}) and the definition of the structure of Abelian category
  on the heart $\mathsf{D}(\theta,\theta+1)$ imply now that the
  morphism $E\rightarrow F$ in (\ref{eq:mutation}) is an epimorphism in
  $\mathsf{D}(\theta,\theta+1)$. This gives an infinite chain of epimorphisms
  which are not isomorphisms, so that the category
  $\mathsf{D}(\theta,\theta+1)$ is indeed not Noetherian.
   
  In order to show that $\mathsf{A}(\theta,\mathsf{P}(\theta)^{-})$ is not
  Noetherian for $\mathsf{P}(\theta)^{-} \ne \emptyset$ we may assume $\theta
  = 0$. If there exists a stable element $\boldsymbol{k}(x) \in
  \mathsf{P}(0)^{-}[1]\subset \mathsf{P}(1)$, where $x\in\boldsymbol{E}$ is a
  smooth point, we have exact sequences
  \begin{equation}
    \label{eq:sequence}
    0 \rightarrow \mathcal{O}(mx)
    \rightarrow \mathcal{O}((m+1)x)
    \rightarrow \boldsymbol{k}(x)
    \rightarrow 0
  \end{equation}
  in $\Coh_{\boldsymbol{E}}$ with arbitrary
  $m\in\mathbb{Z}$.
  Hence the cone of the morphism $\mathcal{O}(mx) \rightarrow
  \mathcal{O}((m+1)x)$ is isomorphic to $\boldsymbol{k}(x)[0]$. Because
  $\boldsymbol{k}(x)[0]$ is an object of $\mathsf{D}^{\le-1}$, with regard to
  the $t$-structure which is defined by $\mathsf{P}({0})^{-}$, we obtain 
  $\tau_{\ge0}(\boldsymbol{k}(x)[0])=0$, which is the cokernel of
  $\mathcal{O}(mx) \rightarrow \mathcal{O}((m+1)x)$ 
  in the Abelian category $\mathsf{A}(0,\mathsf{P}(0)^{-})$, see
  \cite{Asterisque100}, 1.3. 
  Hence, there is an exact sequence
  $$0 \rightarrow \boldsymbol{k}(x)[-1] \rightarrow \mathcal{O}(mx)
  \rightarrow  \mathcal{O}((m+1)x) \rightarrow 0$$
  in $\mathsf{A}(0,\mathsf{P}(0)^{-})$
  and we obtain an infinite chain of epimorphisms 
  $$ \mathcal{O}(x) \rightarrow \mathcal{O}(2x) \rightarrow \mathcal{O}(3x)
  \rightarrow \cdots$$ in the category $\mathsf{A}(0,\mathsf{P}(0)^{-})$,
  which, therefore, is not Noetherian.
  If $\mathsf{P}(0)^{-}[1]$ contains $\boldsymbol{k}(s)$ only, where
  $s\in\boldsymbol{E}$ is the singular point, we proceed as follows. First,
  recall that there exist coherent torsion modules with support at $s$
  which have finite injective dimension, see for example
  \cite{BurbanKreussler}, Section 4. To describe examples of them, we can
  choose a line bundle $\mathcal{L}$ on $\boldsymbol{E}$ and a section
  $\sigma\in H^{0}(\mathcal{L})$, such that the cokernel of
  $\sigma:\mathcal{O}\rightarrow \mathcal{L}$ is a coherent torsion module
  $\mathcal{B}$ of length two with support at $s$. If we embed
  $\boldsymbol{E}$ into $\mathbb{P}^{2}$, such a line bundle $\mathcal{L}$ is
  obtained as the tensor product of the restriction of
  $\mathcal{O}_{\mathbb{P}^{2}}(1)$ with $\mathcal{O}_{\boldsymbol{E}}(-x)$,
  where $x\in\boldsymbol{E}$ is a smooth point. The section $\sigma$
  corresponds to the line in the plane through $x$ and $s$. By twisting with
  $\mathcal{L}^{\otimes m}$ we obtain exact sequences
  $$ 0 \rightarrow \mathcal{L}^{\otimes m} \rightarrow \mathcal{L}^{\otimes
  (m+1)} \rightarrow \mathcal{B} \rightarrow 0$$ in $\Coh_{\boldsymbol{E}}$.
  Because $\mathcal{B}$ is a semi-stable torsion sheaf with support at $s$,
  all its JH-factors are isomorphic to $\boldsymbol{k}(s)$ and we conclude as
  above.
\end{proof}

\begin{proposition}\label{prop:eitheror}
  Let $(\mathsf{D}^{\le0}, \mathsf{D}^{\ge0})$ be a $t$-structure on
  $\Dbcoh(\boldsymbol{E})$ and $B$ a semi-stable indecomposable object
  in $\Dbcoh(\boldsymbol{E})$. Then either $B\in \mathsf{D}^{\le0}$ or $B\in
  \mathsf{D}^{\ge1}$.
\end{proposition}

\begin{proof}
  Let $X\stackrel{f}{\rightarrow} B \stackrel{g}{\rightarrow} Y
  \stackrel{+}{\longrightarrow}$ be a distinguished triangle with $X\in
  \mathsf{D}^{\le0}$ and $Y\in \mathsf{D}^{\ge1}$. Suppose $X\ne 0$ and $Y\ne
  0$ in $\Dbcoh(\boldsymbol{E})$. We decompose both objects into
  indecomposables $X=\bigoplus X_{i}$ and $Y=\bigoplus Y_{j}$. By Lemma
  \ref{lem:tsummands} we have $X_{i}\in \mathsf{D}^{\le0}$ and $Y_{j}\in
  \mathsf{D}^{\ge1}$. If one of the components of the morphisms
  $Y[-1]\rightarrow X=\bigoplus X_{i}$ or $\bigoplus Y_{j}=Y\rightarrow X[1]$
  were zero,  by Lemma \ref{lem:PengXiao} we would obtain a direct summand
  $X_{i}$ or $Y_{j}$ in $B$. Because $B$ was assumed to be indecomposable,
  this implies the claim of the proposition. 

  For the rest of the proof we suppose that all components of these two
  morphisms are non-zero. This implies that $X_{i}$ and $Y_{j}$ are
  non-perfect for all $i,j$. Indeed, if $X_{i}$ were perfect, we could 
  apply Serre duality (\ref{wesPT:i}) to obtain $\Hom(Y,X_{i}[1]) =
  \Hom(X_{i},Y)^{\ast}$, which is zero because $X_{i}\in \mathsf{D}^{\le0}$
  and $Y\in \mathsf{D}^{\ge1}$. The case with perfect $Y_{j}$ can be dealt
  with similarly. 
  Using Lemma \ref{lem:PengXiao} again, it follows that none of the components
  of $f:\bigoplus X_{i} \rightarrow B$ or $g:B\rightarrow \bigoplus Y_{j}$ is
  zero, because none of the $X_{i}$ could be a direct summand of $Y[-1]$ and
  none of the $Y_{j}$ could be a summand of $X[1]$.
  Using Lemma \ref{wesPT:ii}, this implies
  $\varphi_{-}(X_{i}) \le \varphi(B) \le \varphi_{+}(Y_{j})$ for all $i,j$.
  If there exist $i,j$ such that $\varphi_{-}(X_{i}) -
  \varphi_{+}(Y_{j})\not\in \mathbb{Z}$, there exists an integer $k\ge 0$ such
  that $\varphi_{-}(X_{i}[k]) < \varphi_{+}(Y_{j}) < \varphi_{-}(X_{i}[k])
  +1$. Using Proposition \ref{wesPT} (\ref{wesPT:iii}) this implies
  $\Hom(X_{i}[k], Y_{j}) \ne 0$. But, for any 
  integer $k\ge 0$ we have $X_{i}[k]\in \mathsf{D}^{\le0}$ and because
  $Y_{j}\in \mathsf{D}^{\ge1}$, we should have $\Hom(X_{i}[k], Y_{j}) =
  0$. This contradiction implies $\varphi_{-}(X_{i}) -  \varphi_{+}(Y_{j}) \in
  \mathbb{Z}$ for all $i,j$. But, if $k=\varphi_{+}(Y_{j}) -
  \varphi_{-}(X_{i})$, we still have $\Hom(X_{i}[k], Y_{j}) \ne 0$, which
  follows from Proposition \ref{wesPT} (\ref{wesPT:iv}) because $X_{i}$ and
  $Y_{j}$ are not perfect. The conclusion is now that we must have $X=0$ or
  $Y=0$, which implies the claim. 
\end{proof}

\begin{lemma}\label{lem:inequ}
  Let $(\mathsf{D}^{\le0}, \mathsf{D}^{\ge0})$ be a $t$-structure on
  $\Dbcoh(\boldsymbol{E})$. If $F\in \mathsf{D}^{\le0}$ and $G\in
  \mathsf{D}^{\ge1}$, then $\varphi_{-}(F)\ge \varphi_{+}(G)$. 
\end{lemma}

\begin{proof}
  Suppose $\varphi_{-}(F)< \varphi_{+}(G)$. It is sufficient to derive a
  contradiction for indecomposable objects $F$ and $G$.
  Because, for any $k\ge0$, $F[k]\in \mathsf{D}^{\le0}$, we may replace $F$
  by $F[k]$ and can assume $0< \varphi_{+}(G) - \varphi_{-}(F)\le 1$. Now,
  there exists a stable vector bundle $\mathcal{B}$ on $\boldsymbol{E}$ and an
  integer $r$ such that 
  $$\varphi_{-}(F) < \varphi(\mathcal{B}[r]) < \varphi_{+}(G) \le
  \varphi_{-}(F) + 1.$$ 
  By Proposition \ref{prop:eitheror},
  $\mathcal{B}[r]$ is in $\mathsf{D}^{\le0}$ or in
  $\mathsf{D}^{\ge1}$. But, from Proposition \ref{wesPT} (\ref{wesPT:iii}) we
  deduce $\Hom(F, \mathcal{B}[r])\ne 0$ and $\Hom(\mathcal{B}[r],G)\ne 0$. If
  $\mathcal{B}[r]\in \mathsf{D}^{\ge1}$, the first inequality contradicts $F\in
  \mathsf{D}^{\le0}$ and if $\mathcal{B}[r]\in \mathsf{D}^{\le0}$, the second
  one contradicts $G\in \mathsf{D}^{\ge1}$.
\end{proof}

\begin{theorem}\label{thm:tstruc} 
  Let $(\mathsf{D}^{\le0}, \mathsf{D}^{\ge0})$ be a t-structure on
  $\Dbcoh(\boldsymbol{E})$. Then there exists a number $\theta\in \mathbb{R}$
  and a subset 
  $\mathsf{P}(\theta)^{-}\subset \mathsf{P}(\theta)^{s}$, such that  
  $${\sf D}^{\le 0} = \mathsf{D}[\mathsf{P}(\theta)^{-}, \infty)
  \quad\text{ and }\quad
  {\sf D}^{\ge 1} = \mathsf{D}(-\infty,\mathsf{P}(\theta)^{+}].$$
\end{theorem}

\begin{proof}
  From Lemma \ref{lem:inequ} we deduce the existence of $\theta \in
  \mathbb{R}$ such that $\mathsf{D}(\theta,\infty)\subset\mathsf{D}^{\le 0}$
  and $\mathsf{D}(-\infty, \theta)\subset\mathsf{D}^{\ge  1}.$
  If we define
  $\mathsf{P}(\theta)^{-}=\mathsf{P}(\theta)^{s}\cap \mathsf{D}^{\le 0}$ 
  and
  $\mathsf{P}(\theta)^{+}=\mathsf{P}(\theta)^{s}\cap \mathsf{D}^{\ge1}$,
  Proposition \ref{prop:eitheror} implies
  $\mathsf{P}(\theta)^{s} = \mathsf{P}(\theta)^{-}\cup\mathsf{P}(\theta)^{+}$.
  Hence, $\mathsf{D}[\mathsf{P}(\theta)^{-}, \infty)\subset \mathsf{D}^{\le0}$
  and $\mathsf{D}(-\infty,\mathsf{P}(\theta)^{+}]\subset \mathsf{D}^{\ge 1}$. 
  From Proposition \ref{prop:texpl} we know that
  $(\mathsf{D}[\mathsf{P}(\theta)^{-}, \infty),
  \mathsf{D}(-\infty,\mathsf{P}(\theta)^{+}[1]])$ defines a 
  $t$-structure. Now, the statement of the theorem follows.
\end{proof}

\begin{remark}
  In the case of a smooth elliptic curve Theorem \ref{thm:tstruc} was proved in
  \cite{GRK}. If $\theta\not\in\mathsf{Q}$ the heart
  $\mathsf{D}(\theta,\theta+1)$ of the corresponding $t$-structure is a
  finite-dimensional non-Noetherian Abelian category of infinite global
  dimension. In the smooth case, they correspond to the category of holomorphic
  vector bundles on a non-commutative torus in the sense of Polishchuk and
  Schwarz \cite{PolSchw}. It is an interesting problem to find a similar
  interpretation of these Abelian categories in the case of a  singular
  Weierstra{\ss} curve $\boldsymbol{E}$.
\end{remark}

To complete this section we give two applications of Theorem
\ref{thm:tstruc}. The first is a description of the group of exact
auto-equivalences of the triangulated category $\Dbcoh(\boldsymbol{E})$. The
second application is a description of Bridgeland's space of all stability
structures on $\Dbcoh(\boldsymbol{E})$. In both cases, $\boldsymbol{E}$ is an
irreducible curve of arithmetic genus one over $\boldsymbol{k}$.

\begin{corollary}\label{cor:auto}
  There exists an exact sequence of groups
  $$
  \boldsymbol{1} \longrightarrow \Aut^0(\Dbcoh(\boldsymbol{E}))
  \longrightarrow \Aut(\Dbcoh(\boldsymbol{E}))
  \longrightarrow \SL(2,\mathbb{Z}) \longrightarrow \boldsymbol{1}
  $$
  in which $\Aut^0(\Dbcoh(\boldsymbol{E}))$  is generated by tensor products
  with line bundles of degree zero, automorphisms of the curve and the shift
  by $2$. 
\end{corollary}

\begin{proof}
  The homomorphism $\Aut(\Dbcoh(\boldsymbol{E})) \rightarrow \SL(2,\mathbb{Z})$
  is defined by describing the action of an auto-equivalence on
  $\mathsf{K}(\boldsymbol{E})$ in terms of the coordinate functions $(\deg,
  \rk)$.
  That this is indeed in $\SL(2,\mathbb{Z})$ follows, for example, because
  $\Aut(\Dbcoh(\boldsymbol{E}))$ preserves stability and the Euler-form
  \begin{align*}
    \langle \mathcal{F},\mathcal{G}\rangle &=
                       \dim\Hom(\mathcal{F},\mathcal{G}) -
                       \dim\Hom(\mathcal{G},\mathcal{F})\\ 
                       &= \rk(\mathcal{F}) \deg(\mathcal{G}) -
                       \deg(\mathcal{F})\rk(\mathcal{G}) 
  \end{align*}
  for stable and perfect sheaves $\mathcal{F},\mathcal{G}$. 
  Clearly, tensor products with line bundles of degree
  zero, automorphisms of the curve and the shift by $2$ are contained in the
  kernel of this homomorphism. In order to show that the kernel coincides with
  $\Aut^0(\Dbcoh(\boldsymbol{E}))$, we let $\mathbb{G}$ be an arbitrary exact
  auto-equivalence of $\Dbcoh(\boldsymbol{E})$.
  Then, $\mathbb{G}(\Coh_{\boldsymbol{E}})$ is still Noetherian and it is the
  heart 
  of the $t$-structure $(\mathbb{G}(\mathsf{D}^{\le0}),
  \mathbb{G}(\mathsf{D}^{\ge0}))$. 
  From Theorem \ref{thm:tstruc} and Lemma \ref{lem:heart} we know all
  Noetherian hearts of $t$-structures. We obtain
  $\mathbb{G}(\Coh_{\boldsymbol{E}}) = \mathsf{D}(\theta,\theta+1]$ with
  $\mathsf{P}(\theta)\ne \{0\}$. 
  Now, by Corollary \ref{cor:sheaves} there exists
  $\Phi\in\widetilde{\SL}(2,\mathbb{Z})$ which maps
  $\mathsf{D}(\theta,\theta+1]$ to $\mathsf{D}(0, 1]=\Coh_{\boldsymbol{E}}$.
  This implies that the auto-equivalence $\Phi\circ\mathbb{G}$ induces an
  auto-equivalence of the category $\Coh_{\boldsymbol{E}}$. 
  It is well-known that such an auto-equivalence has the form $f^*(\mathcal{L}
  \otimes \ARG)$, where $f:\boldsymbol{E} \rightarrow \boldsymbol{E}$ is an
  isomorphism and $\mathcal{L}$ is a line bundle.
  Note that $f^*(\mathcal{L} \otimes \ARG)$ is sent to the identity in
  $\SL(2,\mathbb{Z})$, if and only if $\mathcal{L}$ is of degree zero.
  The composition of $\Phi\circ\mathbb{G}$ with the inverse of $f^*(\mathcal{L}
  \otimes \ARG)$ satisfies the assumptions of \cite{BondalOrlov}, Prop.~A.3, 
  hence is isomorphic to the identity. 
  Because the kernel of the homomorphism $\widetilde{\SL}(2,\mathbb{Z})
  \rightarrow \SL(2,\mathbb{Z})$, which is induced by the action of
  $\widetilde{\SL}(2,\mathbb{Z})$ on $\Dbcoh(\boldsymbol{E})$ and the above
  homomorphism $\Aut(\Dbcoh(\boldsymbol{E})) \longrightarrow
  \SL(2,\mathbb{Z})$, is generated by the element of
  $\widetilde{\SL}(2,\mathbb{Z})$ which acts as the shift by $2$, the claim
  now follows. 
\end{proof}

For our second application, we recall Bridgeland's definition of stability
condition on a triangulated category \cite{Stability}.

Recall that we set $\mathsf{K}(\boldsymbol{E}) =
\mathsf{K}_{0}(\Coh(\boldsymbol{E})) \cong
\mathsf{K}_{0}(\Dbcoh(\boldsymbol{E}))$. 
Following Bridgeland \cite{Stability}, we call a pair $(W,\mathsf{R})$ a
\emph{stability condition} on $\Dbcoh(\boldsymbol{E})$, if
$$W:\mathsf{K}(\boldsymbol{E})\rightarrow\mathbb{C}$$ 
is a group homomorphism and $\mathsf{R}$ is a compatible slicing of
$\Dbcoh(\boldsymbol{E})$. A \emph{slicing} $\mathsf{R}$ consists of a
collection of full additive subcategories $\mathsf{R}(t) \subset
\Dbcoh(\boldsymbol{E})$, $t\in\mathbb{R}$, such that
\begin{enumerate}
\item $\forall t\in\mathbb{R}\quad \mathsf{R}(t+1) = \mathsf{R}(t)[1]$;
\item If $t_{1}>t_{2}$ and $A_{\nu}\in\mathsf{R}(t_{\nu})$, then
  $\Hom(A_{1},A_{2}) =0$;
\item each non-zero object $X\in\Dbcoh(\boldsymbol{E})$ has a HNF
    \[\xymatrix@C=.4em{
    0\; \ar[rr] && F_{n}X \ar[rr] \ar[dl]_{\cong}&&   F_{n-1}X \ar[rr]
    \ar[dl]&& \dots \ar[rr] &&F_{1}X   \ar[rr]  &&  F_{0}X \ar@{=}[r]\ar[dl] &
    X\\
    & A_{n} \ar[lu]^{+} && A_{n-1} \ar[lu]^{+} &  & & & & &  A_0 \ar[lu]^{+}&
  }\]
  in which $0\ne A_{\nu}\in\mathsf{R}(\varphi_{\nu})$ and
    $\varphi_{n}>\varphi_{n-1}> \ldots > \varphi_{1}>\varphi_{0}$. 
\end{enumerate}
Compatibility means for all non-zero $A\in\mathsf{R}(t)$
$$W(A)\in \mathbb{R}_{>0}\exp(i\pi t).$$
By $\varphi^{\mathsf{R}}$ we denote the phase function on
$\mathsf{R}$-semi-stable objects. Similarly, we denote by
$\varphi^{\mathsf{R}}_{+}(X)$ and 
$\varphi^{\mathsf{R}}_{-}(X)$ the largest, respectively smallest, phase of an
$\mathsf{R}$-HN factor of $X$.

The standard stability condition, which was studied in the previous section,
will always be denoted by $(Z, \mathsf{P})$. This stability condition has a
slicing which is \emph{locally finite}, see \cite{Stability}, Def.\/ 5.7.
A slicing $\mathsf{R}$ is called locally finite, iff there exists $\eta>0$
such that for any $t\in\mathbb{R}$ the quasi-Abelian category
$\mathsf{D}^{\mathsf{R}}(t-\eta, t+\eta)$ is of finite length, i.e. Artinian
and Noetherian.
This category consists of those objects $X\in\Dbcoh(\boldsymbol{E})$ which
satisfy $t-\eta<\varphi^{\mathsf{R}}_{-}(X) \le \varphi^{\mathsf{R}}_{+}(X) <
t+\eta$. 

In order to obtain a good moduli space of stability conditions, Bridgeland
\cite{Stability} requires the stability conditions to be
\emph{numerical}. This means that the central charge $W$ factors 
through the numerical Grothendieck group. This makes sense if for any two
objects $E,F$ of the triangulated category in question, the vector spaces
$\bigoplus_{i} \Hom(E,F[i])$ are finite-dimensional. This condition is not
satisfied for $\Dbcoh(\boldsymbol{E})$, if $\boldsymbol{E}$ is
singular. However, in the case of our interest, we do not need such an extra
condition, because the Grothendieck group $\mathsf{K}(\boldsymbol{E})$ is
sufficiently small. From Lemma \ref{lem:GrothGrp} we know
$\mathsf{K}(\boldsymbol{E}) \cong \mathbb{Z}^{2}$ with generators
$[\mathcal{O}_{\boldsymbol{E}}]$ and $[\boldsymbol{k}(x)]$,
$x\in\boldsymbol{E}$ arbitrary. 
Because $Z(\boldsymbol{k}(x))=-1$ and $Z(\mathcal{O}_{\boldsymbol{E}})=i$, it
is now clear that any homomorphism $W:\mathsf{K}(\boldsymbol{E}) \rightarrow
\mathbb{C}$ can be written as $W(E)=w_{1}\deg(E) + w_{2}\rk(E)$ with $w_{1},
w_{2}\in\mathbb{C}$. Equivalently, if we identify $\mathbb{C}$ with
$\mathbb{R}^{2}$, there exists a $2\times 2$-matrix $A$ such that $W=A\circ Z$.

\begin{definition}
  By $\Stab{\boldsymbol{E}}$ we denote the set of all stability conditions $(W,
  \mathsf{R})$ on $\Dbcoh(\boldsymbol{E})$ for which $\mathsf{R}$ is a locally
  finite slicing. 
\end{definition}

\begin{lemma}\label{lem:notaline}
  For any $(W, \mathsf{R}) \in \Stab(\boldsymbol{E})$ there exists a matrix
  $A\in\GL(2,\mathbb{R})$, such that $W=A\circ Z$.
\end{lemma}

\begin{proof}
  As seen above, there exists a not necessarily invertible matrix $A$ such
  that $W=A\circ Z$. If $A$ were not invertible, there would exist a number
  $t_{0}\in\mathbb{R}$ such that $W(\mathsf{K}(\boldsymbol{E})) \subset
  \mathbb{R}\exp(i\pi t_{0})$. This implies that there may exist a non-zero
  object in $\mathsf{R}(t)$ only if $t-t_{0}\in\mathbb{Z}$. The assumption
  that the slicing $\mathsf{R}$ is locally finite implies now that
  $\mathsf{R}(t)$ is of finite length for any $t\in\mathbb{R}$. On the other
  hand, the heart of the $t$-structure, which is defined by $(W,\mathsf{R})$
  is $\mathsf{R}(t_{0})$ up to a shift. However, in Lemma \ref{lem:heart} we
  determined all Noetherian hearts of $t$-structures on
  $\Dbcoh(\boldsymbol{E})$ and none of them is Artinian. This contradiction
  shows that $A$ is invertible. 
\end{proof}

\begin{lemma}\label{lem:function}
  If $(W,\mathsf{R}) \in \Stab(\boldsymbol{E})$, there exists a unique
  strictly increasing function  $f:\mathbb{R} \rightarrow \mathbb{R}$ with
  $f(t+1) = f(t)+1$ and $\mathsf{R}(t) = \mathsf{P}(f(t))$.
\end{lemma}

\begin{proof}
  By definition, $W(\mathsf{R}(t)) \subset \mathbb{R}_{>0} \exp(i\pi
  t)$. By Lemma \ref{lem:notaline}, there exists a linear isomorphism $A$ such
  that $W=A^{-1}\circ Z$. This implies that there is a function
  $f:\mathbb{R}\rightarrow \mathbb{R}$ such that $Z(\mathsf{R}(t))
  \subset \mathbb{R}_{>0} \exp(i\pi f(t))$. 
  On the other hand, $\mathsf{R}(t)$ is the intersection of two hearts of
  $t$-structures. By Proposition \ref{prop:texpl} these hearts are of the form
  $\mathsf{D}[\mathsf{P}(\theta_{1})^{-}, \mathsf{P}(\theta_{1})^{+}[1]]$ and
  $\mathsf{D}[\mathsf{P}(\theta_{2})^{-}, \mathsf{P}(\theta_{2})^{+}[1]]$ with
  $\theta_{1}\le \theta_{2}$. These have non-empty intersection only if
  $\theta_{2} \le \theta_{1}+1$. Their intersection is contained in
  $\mathsf{D}[\theta_{2},\theta_{1}+1]$, see Figure \ref{fig:intersection}.
  \begin{figure}[hbt]
    \begin{center}
      \setlength{\unitlength}{10mm}
  \begin{picture}(11,5)
    \multiput(0,4)(0.2,0){56}{\line(1,0){0.1}}
    \put(0,1){\line(1,0){11.1}}
    \thicklines
    \put(1.5,2){\line(0,1){2}}\put(1.5,0.8){\makebox(0,0)[t]{$\theta_{2}+1$}}
    \put(1.4,3){\makebox(0,0)[r]{$\mathsf{P}(\theta_{2})^{+}[1]$}}
    \put(5.5,1){\line(0,1){1}}\put(5.5,0.8){\makebox(0,0)[t]{$\theta_{2}$}}
    \put(5.6,1.4){\makebox(0,0)[l]{$\mathsf{P}(\theta_{2})^{-}$}}
    \put(4.5,2.3){\line(0,1){1}}\put(4.5,0.8){\makebox(0,0)[t]{$\theta_{1}+1$}}
    \put(4.4,3){\makebox(0,0)[r]{$\mathsf{P}(\theta_{1})^{+}[1]$}}
    \put(8.5,1){\line(0,1){1.3}}\put(8.5,0.8){\makebox(0,0)[t]{$\theta_{1}$}}
    \put(8.5,3.3){\line(0,1){0.7}}
    \put(8.6,1.5){\makebox(0,0)[l]{$\mathsf{P}(\theta_{1})^{-}$}}
    \thinlines
    \multiput(1.5,1)(0,0.2){5}{\line(0,1){0.1}}
    \multiput(5.5,2)(0,0.2){10}{\line(0,1){0.1}}
    \multiput(4.5,1)(0,0.2){7}{\line(0,1){0.1}}
    \multiput(4.5,3.3)(0,0.2){4}{\line(0,1){0.1}}
    \multiput(8.5,2.3)(0,0.2){5}{\line(0,1){0.1}}
    \put(1.9,4){\line(-2,-3){0.4}}
    \put(2.7,4){\line(-2,-3){1.2}}
    \multiput(1.5,1)(0.8,0){3}{\line(2,3){2}}
    \put(3.9,1){\line(2,3){1.6}}
    \put(4.7,1){\line(2,3){0.8}}
    \put(8.1,4){\line(2,-3){0.4}}
    \put(7.3,4){\line(2,-3){1.2}}
    \multiput(8.5,1)(-0.8,0){3}{\line(-2,3){2}}
    \put(6.1,1){\line(-2,3){1.6}}
    \put(5.3,1){\line(-2,3){0.8}}
    \end{picture}
    \end{center}
    \caption{}\label{fig:intersection}
  \end{figure}

  Moreover, if $\theta_{2}<\theta_{1}+1$, there exist $\alpha,
  \beta\in\mathsf{Q}$ with $\theta_{2}< \alpha < \beta < \theta_{1}+1\le
  \theta_{2}+1$. In this case we have two non-trivial subcategories
  $\mathsf{P}(\alpha)\subset \mathsf{R}(t)$ and $\mathsf{P}(\beta)\subset
  \mathsf{R}(t)$. However, because $0<\beta-\alpha<1$ and 
  $Z(\mathsf{R}(t)) \subset \mathbb{R}_{>0} \exp(i\pi f(t))$, we cannot have
  $Z(\mathsf{P}(\alpha)) \subset \mathbb{R}_{>0} \exp(i\pi\alpha)$ and
  $Z(\mathsf{P}(\beta)) \subset \mathbb{R}_{>0} \exp(i\pi\beta)$.
  Hence, $\theta_{2}=\theta_{1}+1=f(t)$ and we obtain $\mathsf{R}(t) \subset
  \mathsf{P}(f(t))$.  
  From $\mathsf{R}(t+m)=\mathsf{R}(t)[m]$ we easily obtain
  $f(t+m)=f(t)+m$. Moreover, $f(t_{2})=f(t_{1})+m$ with $m\in\mathbb{Z}$
  implies $t_{2}-t_{1}\in\mathbb{Z}$, because the image of $W$ is not
  contained in a line by Lemma \ref{lem:notaline}.

  Next, we show that $f$ is strictly increasing. Suppose $t_{1}<t_{2}$,
  $t_{2}-t_{1}\not\in\mathbb{Z}$ and both $\mathsf{R}(t_{i})$ contain non-zero
  objects $X_{i}$. For any $m\ge0$ we have $\Hom(X_{2}, X_{1}[-m]) = 0$. If
  $f(t_{2}) < f(t_{1})$, we choose $m\ge0$ such that $f(t_{2}) < f(t_{1}) -m <
  f(t_{2}) +1$ and obtain $X_{2}\in\mathsf{P}(f(t_{2}))$ and $X_{1}[-m] \in
  \mathsf{P}(f(t_{1})-m)$. But this implies, by Corollary \ref{cor:sheaves}
  (\ref{cor:iii}), $\Hom(X_{2}, X_{1}[-m]) \ne 0$, a contradiction. Hence, we
  have shown that $f$ is strictly increasing with $f(t+1)=f(t)+1$ and
  $\mathsf{R}(t)\subset \mathsf{P}(f(t))$. In particular, any $\mathsf{R}$-HNF
  is a $\mathsf{P}$-HNF as well. Therefore, all $\mathsf{P}$-semi-stable
  objects are $\mathsf{R}$-semi-stable and we obtain $\mathsf{R}(t) =
  \mathsf{P}(f(t))$. 
\end{proof}

It was shown in \cite{Stability} that the group
$\widetilde{\GL}^{+}(2,\mathbb{R})$ acts naturally on the moduli space of
stability conditions $\Stab(\boldsymbol{E})$. 
This group is the universal cover of $\GL^{+}(2, \mathbb{R})$ and has the
following description:
$$\widetilde{\GL}^{+}(2,\mathbb{R}) = 
\{(A,f) \mid A\in\GL^{+}(2,\mathbb{R}), f:\mathbb{R}\rightarrow \mathbb{R}
\text{ compatible}\},$$
where compatibility means that $f$ is strictly increasing, satisfies
$f(t+1)=f(t)+1$ and induces the same map on $S^{1}\cong\mathbb{R}/2\mathbb{Z}$
as $A$ does on $S^{1}\cong\mathbb{R}^{2}\setminus\{0\}/\mathbb{R}^{\ast}$.
The action is simply $(A,f)\cdot (W,\mathsf{Q})=(A^{-1}\circ W,\mathsf{Q}\circ
f)$. So, this action basically is a relabelling of the slices. 
The following result generalises \cite{Stability}, Thm.\/ 9.1, to the singular
case. 

\begin{proposition}\label{prop:stabmod}
  The action of $\widetilde{\GL}^{+}(2,\mathbb{R})$ on $\Stab(\boldsymbol{E})$
  is simply transitive. 
\end{proposition}

\begin{proof}
  If $(W,\mathsf{R})\in\Stab(\boldsymbol{E})$, the two values
  $W(\mathcal{O}_{\boldsymbol{E}})$ and $W(\boldsymbol{k}(p_{0}))$ determine a
  linear transformation $A^{-1}\in\GL(2, \mathbb{R})$ such that $W=A^{-1}\circ
  Z$, see Lemma \ref{lem:notaline}. 
  By construction, the function $f:\mathbb{R}\rightarrow \mathbb{R}$ of Lemma
  \ref{lem:function} induces the same mapping on
  $S^{1}\cong\mathbb{R}/2\mathbb{Z}$ as $A^{-1}$ does on
  $S^{1}\cong\mathbb{R}^{2}\setminus\{0\}/\mathbb{R}^{\ast}$. Therefore,
  $A\in\GL^{+}(2, \mathbb{R})$ and we obtain $(A,f)\in\widetilde{\GL}^{+}(2,
  \mathbb{R})$ which satisfies $(W,\mathsf{R}) = (A,f)\cdot
  (Z,\mathsf{P})$. Finally, if $(A,f)\cdot (Z,\mathsf{P}) =(Z,\mathsf{P})$ for
  some $(A,f)\in\widetilde{\GL}^{+}(2, \mathbb{R})$, we obtain $f(t)=t$ for
  all $t\in\mathbb{R}$. This implies easily $A=\boldsymbol{1}$. 
\end{proof}

The group $\Aut(\Dbcoh(\boldsymbol{E}))$ acts on $\Stab(\boldsymbol{E})$ by
the rule $$\mathbb{G} \cdot (W, \mathsf{R}) := (\overline{\mathbb{G}}\circ W,
\mathbb{G}(\mathsf{R})).$$ 
Here, $\overline{\mathbb{G}}\in\SL(2,\mathbb{Z})$ is the image of
$\mathbb{G}\in\Aut(\Dbcoh(\boldsymbol{E}))$ under the homomorphism of Corollary
\ref{cor:auto} and $\mathbb{G}(\mathsf{R})(t):= \mathbb{G}(\mathsf{R}(t))$. 
Because automorphisms of $\boldsymbol{E}$ and twists by line 
bundles act trivially on $\Stab(\boldsymbol{E})$, we obtain
$$\Stab(\boldsymbol{E})/\Aut(\Dbcoh(\boldsymbol{E})) \cong \GL^{+}(2,
\mathbb{R})/\SL(2,\mathbb{Z}),$$
which is a $\mathbb{C}^{\times}$-bundle over the coarse moduli space of
elliptic curves. This result coincides with Bridgeland's result in the smooth
case. The main reason for this coincidence seems to be the irreducibility of
the curve. Example \ref{ex:marginalst} below shows that the situation is
significantly more difficult in the case of reducible degenerations of
elliptic curves.  

\begin{remark}\label{rem:common}
  Our results show that singular and smooth Weierstra{\ss} curves
  $\boldsymbol{E}$ share the following properties:
  \begin{enumerate}
  \item A coherent sheaf $\mathcal{F}$ is stable if and only if
    $\End(\mathcal{F}) \cong \boldsymbol{k}$.
  \item Any spherical object is a shift of a stable vector bundle or of a
    structure sheaf $\boldsymbol{k}(x)$ of a smooth point $x\in\boldsymbol{E}$.
  \item The category of semi-stable sheaves of a fixed slope is
    equivalent to the category of coherent torsion sheaves.
    Such an equivalence is induced by an auto-equivalence of
    $\Dbcoh(\boldsymbol{E})$. 
  \item There is an exact sequence of groups\\
    $\boldsymbol{1} \rightarrow \langle \Aut(\boldsymbol{E}),
    \Pic^{0}(\boldsymbol{E}),[2]\rangle \rightarrow
    \Aut(\Dbcoh(\boldsymbol{E})) \rightarrow \SL(2,\mathbb{Z}) \rightarrow
    \boldsymbol{1}.$ 
  \item $\widetilde{\GL}^{+}(2,\mathbb{R})$ acts transitively on
    $\Stab(\boldsymbol{E})$. 
  \item $\Stab(\boldsymbol{E})/\Aut(\Dbcoh(\boldsymbol{E})) \cong \GL^{+}(2,
  \mathbb{R})/\SL(2,\mathbb{Z}).$ 
  \end{enumerate}
\end{remark}

\begin{example}\label{ex:marginalst}
  Let $C_{2}$ denote a reducible curve which has two components, both
  isomorphic to $\mathbb{P}^{1}$ and which intersect transversally at two 
  distinct points. This curve has arithmetic genus one and appears as a
  degeneration of a smooth elliptic curve.
  On this curve, there exists a line bundle $\mathcal{L}$ which
  fails to be stable with respect to some stability conditions. To construct
  an explicit example, denote by $\pi:\widetilde{C}_{2}\rightarrow C_{2}$ the
  normalisation, so that $\widetilde{C}_{2}$ is the disjoint union of two
  copies of $\mathbb{P}^{1}$. There is a $\boldsymbol{k}^{\times}$-family of
  line bundles whose pull-back to $\widetilde{C}_{2}$ is
  $\mathcal{O}_{\mathbb{P}^{1}}$ on one component and
  $\mathcal{O}_{\mathbb{P}^{1}}(2)$ on the other. The element in
  $\boldsymbol{k}^{\times}$ 
  corresponds to a gluing parameter over one of the two singularities. Let
  $\mathcal{L}$ denote one such line bundle.
  If $i_{\nu}:\mathbb{P}^{1}\rightarrow C_{2},\;\nu=1,2$ denote the embeddings
  of the two components, we fix notation so that $i_{1}^{\ast}\mathcal{L}
  \cong \mathcal{O}_{\mathbb{P}^{1}}$ and $i_{2}^{\ast}\mathcal{L} \cong
  \mathcal{O}_{\mathbb{P}^{1}}(2)$. There is an exact sequence of coherent
  sheaves on $C_{2}$ 
  \begin{equation}\label{eq:linebundle}
  0\rightarrow i_{2\ast} \mathcal{O}_{\mathbb{P}^{1}} \rightarrow \mathcal{L}
  \rightarrow i_{1\ast} \mathcal{O}_{\mathbb{P}^{1}} \rightarrow 0.
  \end{equation}
  Moreover, the only non-trivial quotients of $\mathcal{L}$ are
  $\mathcal{L}\twoheadrightarrow i_{1\ast} \mathcal{O}_{\mathbb{P}^{1}}$
  and  
  $\mathcal{L}\twoheadrightarrow i_{2\ast} \mathcal{O}_{\mathbb{P}^{1}}(2)$.

  For arbitrary positive real numbers $a,b$ we define a centred
  slope-function $W_{a,b}$ on the category $\Coh_{C_{2}}$ by 
  $$W_{a,b}(F):= -\deg(F) + i(a\cdot \rk(i_{1}^{\ast}F) + b\cdot
  \rk(i_{2}^{\ast} F)),$$ 
  where $\deg(F)=h^{0}(F) - h^{1}(F)$. For example,
  \begin{align*}
    W_{a,b}(i_{1\ast}\mathcal{O}_{\mathbb{P}^{1}}(d)) &=  -d-1+ia
    \quad\text{ and }\\
    W_{a,b}(i_{2\ast}\mathcal{O}_{\mathbb{P}^{1}}(d)) &=  -d-1+ib.
  \end{align*}
  Using the exact sequence (\ref{eq:linebundle}), we obtain
  $W_{a,b}(\mathcal{L}) = -2+i(a+b)$.
  Using results of \cite{Rudakov}, it is easy to see that $W_{a,b}$ has the
  Harder-Narasimhan property in the sense of \cite{Stability}. Hence, by
  \cite{Stability}, Prop.\/ 5.3, $W_{a,b}$ defines a stability condition on
  $\Dbcoh(C_{2})$. With respect to this stability condition, the line bundle
  $\mathcal{L}$ is stable precisely when $2/(a+b) < 1/a$, which is
  equivalent to $a<b$. It is semi-stable, but not stable, if $b=a$. If $a>b$,
  $\mathcal{L}$ is not even semi-stable. 
\end{example}

This example illustrates an interesting effect, which was not available on an
irreducible curve of arithmetic genus one. It is an interesting problem to
describe the subset in $\Stab(\boldsymbol{E})$ for which a given line
bundle $\mathcal{L}$ is semi-stable. This is a closed subset, see
\cite{Stability}. Physicists call the boundary of this set the line of
marginal stability, see e.g. \cite{AspinwallDouglas}. The example above
describes the intersection of this set with a two-parameter family of
stability conditions in $\Stab(\boldsymbol{E})$.

\begin{remark}
  In the case of an irreducible curve of arithmetic genus one, we have shown
  in Proposition \ref{prop:spherical} that $\Aut(\Dbcoh(\boldsymbol{E}))$ acts
  transitively on the set of all spherical objects on
  $\boldsymbol{E}$. Polishchuk \cite{YangBaxter} conjectured that this is
  likewise true in the case of reducible curves with trivial dualising sheaf.
  However, on the curve $C_{2}$ there exists a
  spherical complex which has non-zero cohomology in two different degrees, see
  \cite{BuBu}. This indicates that the reducible case is more difficult and
  involves new features.
\end{remark}

\end{document}